\magnification=\magstep1
\hsize=16truecm
 
\input amstex
\TagsOnRight
\parindent=20pt
\parskip=2.5pt plus 1.2pt
\define\({\left(}
\define\){\right)}
\define\[{\left[}
\define\]{\right]}
\define\e{\varepsilon}

\define\const{\text{\rm const.}\,}

\define\supp {\sup\limits}

\define\summ{\sum\limits}
\define\prodd{\prod\limits}
\define\limm{\lim\limits}

\define\Var{\text{\rm Var}\,}

\font\kisit=cmti8
 
\define\sumn{\operatornamewithlimits{ {\sum}^{\prime(\text{\kisit n})}}
}
\define\sumnl{\operatornamewithlimits{ {\sum}^{\prime(\text{\kisit
n,$\,$ L })}} }
 
\hbox{}
\vskip2truecm
 
\centerline{\bf ON A MULTIVARIATE VERSION OF}
\centerline {\bf BERNSTEIN'S INEQUALITY}
\medskip
\centerline{\it P\'eter Major}
\centerline{Alfr\'ed R\'enyi Mathematical Institute of the Hungarian
Academy of Sciences}
\centerline{Budapest, P.O.B. 127 H--1364, Hungary, e-mail:
major\@renyi.hu}
\medskip
 
{\narrower \noindent {\it Summary:}\/ We prove a multivariate
version of Bernstein's inequality about the probability that
degenerate $U$-statistics take a value larger than some number~$u$.
This is an improvement of former estimates for the same problem
which yields an asymptotically sharp estimate for not too large
numbers~$u$. This paper also contains an analogous bound about the
distribution of multiple Wiener-It\^o integrals. Their comparison
shows that our results are sharp. The proofs are based on good
estimates about high moments of multiple random integrals. They are
obtained by means of a diagram formula which enables us to express
the product of multiple random integrals as the sum of such
expressions. \par}
 
\vfill\eject

\beginsection 1. Introduction.
 
Let us consider a sequence of iid.\ random variables
$\xi_1,\xi_2,\dots$, on a measurable space $(X,\Cal X)$ with
some distribution $\mu$ together with a real valued function
$f=f(x_1,\dots,x_k)$  of $k$ variables defined on the $k$-th
power $(X^k,\Cal X^k)$ of the space $(X,\Cal X)$ and  define
with their help the $U$-statistics~$I_{n,k}(f)$, $n=k,k+1,\dots$,
$$
I_{n,k}(f)=\frac1{k!}\summ\Sb 1\le j_s\le n,\; s=1,\dots, k\\
j_s\neq j_{s'} \text{ if } s\neq s'\endSb
f\(\xi_{j_1},\dots,\xi_{j_k}\).   \tag1.1
$$
We are interested in good estimates on probabilities of the type
$P\(n^{-k/2}k!|I_{n,k}(f)|>u\)$ under appropriate conditions.
 
Arcones and Gin\'e in [1] have proved an inequality which can be
written in a slightly different but equivalent form as
$$
\aligned
P\(k!n^{-k/2}|I_{n,k}(f)|>u\)\le
c_1\exp&\left\{-\dfrac{c_2u^{2/k}}{\sigma^{2/k}\(1+c_3
\(un^{-k/2}\sigma^{-(k+1)}\)^{2/k(k+1)}\)} \right\} \\
&\qquad \text{for all }u>0
\endaligned\tag1.2
$$
with some universal constants $c_1$, $c_2$ and $c_3$ depending only
on the order $k$ of the $U$-statistic $I_{n,k}(f)$ defined in (1.1)
if the function $f$ satisfies the conditions
$$
\align
\|f\|_\infty&=\sup_{x_j\in X, \,1\le j\le k}|f(x_1,\dots,x_k)|\le 1,
\tag1.3 \\
\|f\|^2_2&=\int f^2(x_1,\dots,x_k)\mu(\,dx_1)\dots\mu(\,dx_k)
\le\sigma^2, \tag1.4
\endalign
$$
and it is canonical with respect to the probability measure $\mu$, i.e.
$$
\align
\int f(x_1,\dots,x_{j-1},u,x_{j+1},\dots,x_k)\mu(\,du)=0\quad
&\text{for all } 1\le j\le k\\
&\quad \text{and \ } x_s\in X, \; s\in\{1,\dots k\}\setminus \{j\}.
\endalign
$$
A $U$-statistic defined in (1.1) with the help of a canonical
function $f$ is called degenerate in the literature. A degenerate
$U$-statistic is the natural multivariate version of sums of iid.
random variables with expectation zero.
 
Arcones and Gin\'e called their estimate (1.2) a new Bernstein-type
inequality. The reason for such a name is that the original Bernstein
inequality (see e.g.~[3], 1.3.2 Bernstein inequality) states relation
(1.2) in the case $k=1$ with the constants $c_1=2$, $c_2=\frac12$
and $c_3=\frac13$ if the function $f(x)$ satisfies the conditions
$\sup|f(x)|\le1$, $\int f(x)\mu(\,dx)=0$ and $\int f^2(x)\mu(\,dx)
\le\sigma^2$.
 
Let us fix a number $C>0$. Formula (1.2) states in particular
that for all numbers $0\le u<Cn^{k/2}\sigma^{k+1}$ and degenerate
$U$-statistics $I_{n,k}(f)$ of order $k$ with a kernel function
$f$ satisfying relations (1.3) and (1.4) the inequality
$P\(n^{-k/2}k!|I_{n,k}(f)|>u\)\le A\exp\left\{-B\(\frac
u\sigma\)^{2/k}\right\}$ holds with some appropriate constants
$A=A(C,k)$ and $B=B(C,k)$ depending only on the fixed number $C$ and
the order $k$ of the degenerate $U$-statistics. This inequality can
be interpreted in the following way: Let us take a random variable
$\eta$ with standard normal distribution. Then
$P\(n^{-k/2}k!|I_{n,k}(f)|>u\)\le A
P\(\sigma\left|\frac\eta{\sqrt{2B}}\right|^k>u\)$, at least for
$0\le u< Cn^{k/2}\sigma^{k+1}$. Let us also observe that under
condition~(1.4) the variance of $n^{-k/2}k!I_{n,k}(f)$ is bounded
by $k!\sigma^2$, and if the kernel function $f$ is symmetric
and there is identity instead of inequality in (1.4), then
$\limm_{n\to\infty}\Var\(n^{-k/2}k!I_{n,k}(f)\)=k!\sigma^2$.
 
In the above discussion we have considered the probability
$P\(n^{-k/2}k!|I_{n,k}(f)|>u\)$ only for
$0\le u<Cn^{k/2}\sigma^{k+1}$, while formula (1.2) yields an
estimate for such a probability for all $u>0$. On the other hand,
as I shall show later, the above restriction of the parameter~$u$
does not mean an important loss of information.
 
Bernstein's inequality yields an analogous estimate in the case
of degenerate $U$-statistics of order 1, i.e. for sums
$\summ_{j=1}^nf(\xi_j)$ with a sequence of iid.\ random variables
$\xi_1,\dots,\xi_n$ and a function $f(x)$ whose absolute value is
bounded by 1, $Ef(\xi_1)=0$ and $Ef(\xi_1)^2=\sigma^2$. (Actually,
Bernstein's inequality is more general, it also yields a bound for
the distribution of a sum of independent, not necessarily identically
distributed random variables.) But Bernstein's inequality also
contains some additional information. It states that if $0\le u<\e
n^{1/2}\sigma^2$ with a small $\e>0$, then
$P\(n^{-1/2}\left|\summ_{j=1}^nf(\xi_j)\right|>u\)\le
P((1-C\e)\sigma|\eta|>u)$ with an appropriate constant $C>0$.
Since $n^{-1/2}\summ_{j=1}^nf(\xi_j)$ has expectation zero and
variance $\sigma^2$ the above inequality can be interpreted in such
a way that at not too large values~$u$ the distribution function of
the normalized sum $n^{-1/2}\summ_{j=1}^nf(\xi_j)$ can be bounded by
the distribution of a normal random variable with expectation zero
and only slightly smaller variance. The main goal of this paper
is to show that a similar estimate holds for degenerate
$U$-statistics of any order.
 
To carry out such a program first we have to find a good
multivariate analog of Gaussian random variables. It is natural to
consider multiple Wiener--It\^o integrals which also appear as the
limit of normalized degenerate $U$-statistics as the sample size
tends to infinity. (See e.g.~[4]). We shall prove an estimate about
the distribution of multiple Wiener--It\^o integrals in Theorem~1
and show in Example~2 that this estimate is sharp. The main result
of this paper is Theorem~3 which yields an estimate about the tail
behaviour of degenerate $U$-statistics. Its comparison with
Theorem~1 shows that Theorem~3 provides an asymptotically sharp
estimate on the tail distribution of a degenerate $U$-statistic
for not too large values.
 
To formulate Theorem~1 let us take a $\sigma$-finite measure $\mu$
on the space $(X,\Cal X)$ and a white noise $\mu_W$ with counting
measure $\mu$ on $(X,\Cal X)$, i.e.\ a set of jointly Gaussian
random variables $\mu_W(A)$, $A\in\Cal X$, such that $E\mu_W(A)=0$,
$E\mu_W(A)\mu_W(B)=\mu(A\cap B)$ for all $A\in\Cal X$ and $B\in\Cal X$.
(We also need the identity $\mu_W(A\cup B)=\mu_W(A)+\mu_W(B)$ with
probability~1 if $A\cap B=\emptyset$, but this is the consequence of
the previous properties of the white noise. Indeed, they imply that
$E\[\mu_W(A\cup B)-(\mu_W(A)+\mu_W(B))\]^2=0$ if $A\cap B=\emptyset$,
hence the desired identity holds.)
The $k$-fold Wiener--It\^o integral of a function $f$
$$
J_{\mu,k}(f)=\frac1{k!}\int
f(x_1,\dots,x_k)\mu_W(\,dx_1)\dots\mu_W(dx_k) \tag1.5
$$
can be defined with respect to a white noise $\mu_W$ with counting
measure $\mu$ if $f$ is a measurable function on the space
$(X^k,\Cal X^k)$, and it satisfies relation (1.4) with some
$\sigma^2<\infty$. (See e.g~[6] or~[7].) The expression $J_{\mu,k}(f)$
in formula (1.5) will be called a Wiener--It\^o integral of order $k$.
Our first result is the following estimate which is an improvement of
the upper bound given in Theorem~6.6 of~[7].
 
\medskip\noindent
{\bf Theorem 1.} {\it Let us consider a $\sigma$-finite measure $\mu$
on a measurable space together with a white noise $\mu_W$
with counting measure $\mu$. Let us have a real-valued function
$f(x_1,\dots,x_k)$ on the space $(X^k,\Cal X^k)$ which satisfies
relation (1.4) with some $\sigma^2<\infty$. Take the random
integral $J_{\mu,k}(f)$ introduced in formula (1.5). This random
integral satisfies the inequality
$$
P(k!|J_{\mu,k}(f)|>u)\le C \exp\left\{-\frac12\(\frac
u\sigma\)^{2/k}\right\}\quad \text{for all } u>0 \tag1.6
$$
with an appropriate constant $C=C(k)>0$ depending only on the
multiplicity $k$ of the integral.}
\medskip
 
The following example shows that the estimate of Theorem~1 is sharp.
\medskip\noindent
{\bf Example 2.} {\it Let us have a $\sigma$-finite measure $\mu$ on
some measure space $(X,\Cal X)$ together with a white noise $\mu_W$
on $(X,\Cal X)$ with counting measure~$\mu$. Let $f_0(x)$ be a real
valued function on $(X,\Cal X)$ such that $\int f_0(x)^2\mu(\,dx)=1$,
and take the function $f(x_1,\dots,x_k)=\sigma f_0(x_1)\cdots f_0(x_k)$
with some number $\sigma>0$ and the Wiener--It\^o integral
$J_{\mu,k}(f)$ introduced in formula (1.5).
 
Then the relation
$\int f(x_1,\dots,x_k)^2\,\mu(\,dx_1)\dots\,\mu(\,dx_k)=\sigma^2$
holds, and the random integral $J_{\mu,k}(f)$ satisfies the inequality
$$
P(k!|J_{\mu,k}(f)|>u)\ge \frac{\bar C}{\(\frac u\sigma\)^{1/k}+1}
\exp\left\{-\frac12\(\frac
u\sigma\)^{2/k}\right\}\quad \text{for all } u>0 \tag1.7
$$
with some constant $\bar C>0$.}
\medskip\noindent
{\it Proof of the statement of Example 2:}\/ We may restrict our
attention to the case $k\ge2$. It\^o's formula (see~[6] or~[7]) states
that the random variable $k!J_{\mu,k}(f)$ can be expressed as
$k!J_{\mu,k}(f)=\sigma H_k\(\int f_0(x)\mu_W(\,dx)\)=\sigma H_k(\eta)$,
where $H_k(x)$ is the $k$-th Hermite polynomial with leading
coefficient~1, and $\eta=\int f_0(x)\mu_W(\,dx)$ is a standard normal
random variable. Hence we get by exploiting that the coefficient of
$x^{k-1}$ in the polynomial $H_k(x)$ is zero that
$P(k!|J_{\mu,k}(f)|>u)=P(|H_k(\eta)| \ge\frac u\sigma)\ge
P\(|\eta^k|-D|\eta^{k-2}|>\frac u\sigma\)$ with a sufficiently large
constant $D>0$ if $\frac u\sigma>1$. There exist such positive
constants $A$ and $B$ that $P\(|\eta^k|-D|\eta^{k-2}|>\frac u\sigma\)
\ge P\(|\eta^k|>\frac u\sigma+A\(\frac u\sigma\)^{(k-2)/k}\)$ if
$\frac u\sigma>B$.
 
Hence
$$
P(k!|J_{\mu,k}(f)|>u)\ge
P\(|\eta|>\(\frac u\sigma\)^{1/k}\(1+A\(\frac u\sigma\)^{-2/k}\)\)
\ge  \frac{\bar C \exp\left\{-\frac12\(\frac u\sigma\)^{2/k}\right\}}
{\(\frac u\sigma\)^{1/k}+1}
$$
with an appropriate $\bar C>0$ if $\frac u\sigma>B$. Since
$P(k!|J_{\mu,k}(f)|>0)>0$, the above inequality also holds
for $0\le \frac u\sigma\le B$ if the constant $\bar C>0$ is chosen
sufficiently small. This means that relation (1.7) holds.
 
The main result of this paper is the following
\medskip\noindent
{\bf Theorem 3.} {\it Let $\xi_1,\dots,\xi_n$ be a sequence of iid.
random variables on a space $(X,\Cal X)$ with some distribution~$\mu$.
Let us consider a function $f(x_1,\dots,x_k)$  canonical with
respect to the measure~$\mu$ on the space $(X^k,\Cal X^k)$ which
satisfies conditions (1.3) and (1.4) with some $0<\sigma^2\le1$
together with the degenerate $U$-statistic $I_{n,k}(f)$ with this
kernel function. There exist some constants $A=A(k)>0$ and
$B=B(k)>0$ depending only on the order $k$ of the $U$-statistic
$I_{n,k}(f)$ such that
$$
P(k!n^{-k/2}|I_{n,k}(f)|>u)\le A\exp\left\{-\frac{u^{2/k}}{2\sigma^{2/k}
\(1+B\(un^{-k/2}\sigma^{-(k+1)}\)^{1/k}\)}\right\} \tag1.8
$$
for all $0\le u\le n^{k/2}\sigma^{k+1}$.}
\medskip\noindent
{\it Remark:} Actually, the universal constant $B>0$  can be chosen
independently of the order $k$ of the degenerate $U$-statistic
$I_{n,k}(f)$ in inequality (1.8).
\medskip
Theorem 3 states in particular that if $0<u\le\e n^{k/2}\sigma^{k+1}$
with a sufficiently small $\e>0$, then $P(k!n^{-k/2}|I_{n,k}(f)|>u)\le
A\exp\left\{-\frac{1-C\e^{1/k}}2\(\frac u\sigma\)^{2/k}\right\}$
with some universal constants $A>0$ and $C>0$ depending only on the
order~$k$ of the $U$-statistic $I_{n,k}(f)$. A  comparison of this
result with Theorem~1 and Example~2 shows that for small $\e>0$ this
estimate yields the right order in the exponent in first order. This
means that for not too large numbers $u$ inequality (1.8) yields an
asymptotically optimal estimate.
 
To understand the previous statement better we can make the
following observation: Let us have a probability measure $\mu$ on
some measurable space $(X,\Cal X)$ together with a sequence of iid.\
random variables $\xi_1,\xi_2,\dots$, with distribution $\mu$, a
real-valued function $f_0(x)$ on $(X,\Cal X)$ such that $\int
f_0(x)\mu(\,dx)=0$, $\int f_0^2(x)\mu(\,dx)=1$ and a real
number~$\sigma$. Let us introduce the function
$f(x_1,\dots,x_k)=\sigma f_0(x_1)\cdots f_0(x_k)$ on $(X^k,\Cal X^k)$
and the $U$-statistics $I_{n,k}(f)$, $n=1,2,\dots$, of order~$k$
defined in formula~(1.1) with this function~$f$. Then the
$U$-statistics $I_{n,k}(f)$ are degenerate, the normalized
$U$-statistics $n^{-k/2}I_{n,k}(f)$ converge in distribution to
the multiple Wiener--It\^o integral $J_{\mu,k}(f)$ introduced in
Example~2 (with the same measure $\mu$ and function $f$ which
appears in the definition of $I_{n,k}(f)$) as $n\to\infty$, (see
e.g.~[4]), and this Wiener--It\^o integral satisfies relation~(1.7).
If the supremum of the function~$f$ is bounded by~1, then Theorem~3
can be applied for the $U$-statistics $k!n^{-k/2}I_{n,k}(f)$, and
the above considerations indicate that for not too large values~$u$
the estimate~(1.8) is sharp.
 
Our goal was to find such an estimate about the distribution function
of degenerate $U$-statistics which is asymptotically optimal for not
too large values of this function. Inequality (1.8) has this property.
In this respect it is similar to Bernstein's inequality which yields
such an estimate in the special case $k=1$. On the other hand,
inequality~(1.2) proved in~[1], does not supply such a bound.
Moreover, the method of paper~[1] seems not to be strong enough to
yield such an estimate. I return to this question later.
 
Let us remark that relation (1.2) yields a bound for the
tail-distribution of a degenerate $U$-statistic for all
numbers~$u>0$, while formula~(1.8) holds under the condition
$0\le u\le n^{k/2}\sigma^{k+1}$. Nevertheless, formula~(1.8)
implies an estimate also for $u>n^{k/2}\sigma^{k+1}$ which is
not weaker than the estimate~(1.2) (at least if we do not fix the
universal constants in these estimates). To see this let us first
observe that in the case $n^{k/2}\ge u>n^{k/2}\sigma^{k+1}$ relation
(1.8) holds with $\bar\sigma=\(u{n^{-k/2}}\)^{1/{(k+1)}}$
and it yields that
$P(k!n^{-k/2}|I_{n,k}(f)|>u)\le A\exp\left\{-\frac1{2(1+B)^{1/k}}
\(\frac u{\bar\sigma}\)^{2/k}\right\}
=Ae^{-(u^2n)^{1/(k+1)}/2(1+B)^{1/k}}$. On the other hand,
$\sigma^{2/k}\(1+c_3\(un^{-k/2}\sigma^{-(k+1)}\)^{2/k(k+1)}\)\ge
c_3u^{2k/(k+1)}n^{-1/k}$, hence the right-hand side of (1.2)
can be bounded from below by $c_1e^{-c_2(u^2n)^{1/(k+1)}/c_3}$.
Thus relation (1.8) implies relation (1.2) if $n^{k/2}\ge
u>n^{k/2}\sigma^{k+1}$ with possibly worse constants $\bar c_1=A$,
$c_2$ and $\bar c_3=2c_2(1+B)^{1/k}$. If $u>n^{k/2}$, then
the left-hand side of (1.2) equals zero because of the boundedness
of the function~$f$, and relation (1.2) clearly holds.
 
Theorem 3 shows some analogy with large deviation results about
the average of iid. random variables. If we fix some number larger
than the expected value of the average of some iid. random variables,
then by the large deviation theory this average can be larger than
this number only with exponentially small probability.
The term in the exponent of the formula expressing
this probability strongly depends on the distribution of the random
variables whose average is taken. But if the above probability is
considered at a level only slightly greater than the expectation of
the average, then this term in the exponent can be well approximated
by the value suggested by the central limit theorem. A similar result
holds for the distribution of normalized degenerate $U$-statistics,
$n^{-k/2}k!I_{n,k}(f)$. In the case $0\le u\le \const
n^{k/2}\sigma^{k+1}$, with $\sigma^2=Ef^2(\xi_1,\dots,\xi_k)$ we can
get a large deviation type estimate for the probability
$P(n^{-k/2}k!I_{n,k}(f)>u)$. If $0\le u\le \e n^{k/2}\sigma^{k+1}$
with a small $\e>0$, then we can say more. In this case such an
estimate can be given which is suggested by the behaviour of
appropriate non-linear functionals of Gaussian processes.
 
Let me also remark that in the case $u\gg n^{k/2}\sigma^{k+1}$
formula (1.8) (or (1.2)) yields only a rather weak estimate for
the probability $P(n^{-k/2}k!I_{n,k}(f)>u)$ for a degenerate
$U$-statistic of order~$k$ with a kernel function~$f$ satisfying
relation~(1.3) and~(1.4). The weakness of our estimate in this case
has a deeper cause. In Examples~3.3 and~8.6 of the work~[10] I have
presented such examples for degenerate $U$-statistics of order~1
or~2 with a kernel function~$f$ satisfying relations~(1.3) and~(1.4)
for which the lower bounds
$P(n^{-1/2}I_{n,1}(f)>u)\ge\exp\left\{-An^{1/2}u\log\frac
u{n\sigma^2}\right\}$, and $P(n^{-1}I_{n,2}(f)>u)\ge
\exp\left\{-An^{1/3}u^{2/3}\log\(\frac u{n\sigma^3}\)\right\}$ hold
if $B_2n^{1/2}\ge u\ge B_1n^{1/2}\sigma^2$ or $B_2n\ge u\ge
B_1n\sigma^3$ respectively with a sufficiently large $B_1>0$ and some
appropriate $0<B_2<1$. Similar examples of degenerate
$U$-statistics could also be constructed for any order~$k$. Thus
there are such degenerate $U$-statistics of order $k$ with a kernel
function satisfying relation (1.3) and (1.4) with some $\sigma>0$,
whose tail distribution have an essentially different
behaviour for $u<n^{k/2}\sigma^{k+1}$ and $u\gg n^{k/2}\sigma^{k+1}$.
 
There is another sort of interesting generalization of Bernstein's
inequality. I would refer to a recent work of C.~Houdr\'e and
P.~Reynaud--Bouret~[5], where good estimates are given for the
distribution of degenerate $U$-statistics of order~2, but in that
paper a more general model is considered. It deals with a
natural object we can call generalized $U$-statistic in the special
case $k=2$. Generalized $U$-statistics can be defined similarly to
classical $U$-statistics with the difference that the underlying
independent random variables $\xi_1,\xi_2,\dots$ may be not
identically distributed, and the terms in the sum (1.1) are of the
form $f_{j_1,\dots,j_k}(\xi_{j_1},\dots,\xi_{j_k})$. If the functions
$f_{j_1,\dots,j_k}$ are canonical, then we can speak about generalized
degenerate $U$-statistics. (The notion of canonical functions can be
generalized to this case in a natural way.) The problem about the
distribution of generalized degenerate $U$-statistics can be
considered as a multivariate version  of the problem about the
distribution of sums of independent, but not necessarily identically
distributed random variables with expectation zero. Here we do not
discuss this question, although it is very interesting. The most
essential part of this problem seems to be to find the right
formulation of the estimate we have to prove. A good estimate on the
distribution of generalized degenerate $U$-statistics has to depend
beside the variance of the $U$-statistics on different quantities
which still should be found.
 
It is natural to expect that generalized degenerate $U$-statistics
$I_{n,k}(f)$ of order~$k$ (without normalization) satisfy the
inequality
$$
P(|I_{n,k}(f)|>u)<A\exp\left\{-C\(\frac u{V_n}\)^{2/k}\right\} \tag1.9
$$
with some universal constants $A=A(k)>0$ and $C=C(k)>0$ in a
relatively large interval for the parameter~$u$, where $V^2_n$
denotes the variance of $I_{n,k}(f)$. An essential
problem is to find a relatively good constant $C$ and to determine
the interval $0<u<D_n$, where the estimate~(1.9) holds. The result
of this paper states that in the case of classical degenerate
$U$-statistics (1.9) holds in the interval $[0,D_n]$ with $D_n=\const
n^k\sigma^{k+1}$, where $\sigma^2=Ef(\xi_1,\dots,\xi_k)^2$.
For $k=1$ this means that relation (1.9) holds in the interval $0\le
u\le V_n^2$. But it is not clear what corresponds in the general case
to the right end-point $D_n=\const n^k\sigma^{k+1}$ of the interval
where the estimate (1.9) should hold.
 
This paper consists of six sections and an Appendix. In Section~2
the method of the proofs is explained. Our results will be proved
by means of a good estimate on high (but not too high) moments
of the random variables we are investigating. These estimates are
obtained by means of a diagram formula which enables us to express
product of stochastic integrals or degenerate $U$-statistics as
a sum of such expressions. Section~3 contains the proof of
Theorem~1. We formulate a version of the diagram formula about
the product of two degenerate $U$-statistics in Section~4. In
Section~5 this result will be generalized to the product of $L\ge2$
degenerate $U$-statistics, and an estimate is given about the
$L_2$-norm of the kernel functions appearing in the $U$-statistics
of this result. Theorem~3 will be proved in Section~6. The diagram
formula about the product of two degenerate $U$-statistics is proved
in the Appendix.
 
\beginsection 2. The idea of the proof.
 
Theorem 1 will be proved by means of the following
\medskip\noindent
{\bf Proposition A.} {\it Let the conditions of Theorem~1 be satisfied
for a multiple Wiener--It\^o integral $J_{\mu,k}(f)$ of order~$k$.
Then, with the notations of Theorem~1, the inequality
$$
E\(k!|J_{\mu,k}(f)|\)^{2M}\le 1\cdot3\cdot5\cdots
(2kM-1)\sigma^{2M}\quad\text {for all }M=1,2,\dots       \tag2.1
$$
holds.}\medskip
By the Stirling formula Proposition~A implies that
$$
E(k!|J_{\mu,k}(f)|)^{2M}\le \frac{(2kM)!}{2^{kM}(kM)!}\sigma^{2M}\le
A\(\frac2e\)^{kM}(kM)^{kM}\sigma^{2M} \tag2.2
$$
for any $A>\sqrt2$ if $M\ge M_0=M_0(A)$.
The following Proposition~B which will be applied in the proof of
Theorem~3 states a similar, but weaker inequality for the moments
of normalized degenerate $U$-statistics.
\medskip\noindent
{\bf Proposition B.} {\it Let us consider a degenerate $U$-statistic
$I_{n,k}(f)$ of order $k$ with sample size $n$ and with a kernel
function $f$ satisfying relations (1.3) and (1.4) with some
$0<\sigma^2\le1$. Fix a positive number $\eta>0$. There exist some
universal constants $A=A(k)>\sqrt2$, $C=C(k)>0$ and $M_0=M_0(k)\ge1$
depending only on the order of the $U$-statistic $I_{n,k}(f)$ such that
$$
\aligned
E\(n^{-k/2}k!I_{n,k}(f)\)^{2M}&\le A\(1+C\sqrt\eta\)^{2kM}
\(\frac2e\)^{kM}\(kM\)^{kM}\sigma^{2M} \\
&\qquad \text{for all integers } M \text{ such that } kM_0\le kM\le
\eta n\sigma^2.
\endaligned \tag2.3
$$
The constant $C=C(k)$ in formula (2.3) can be chosen e.g. as
$C=2\sqrt2$ which does not depend on the order $k$ of the
$U$-statistic $I_{n,k}(f)$.}\medskip
 
Let us remark that formula (2.1) can be reformulated as
$E(k!|J_{\mu,k}(f)|)^{2M}\le E(\sigma\eta^k)^{2M}$, where $\eta$
is a standard normal random variable. Theorem~1 states that
the tail distribution of $k!|J_{\mu,k}(f)|$ satisfies an estimate
similar to that of $\sigma|\eta|^k$. This will simply
follow from Proposition~A and the Markov inequality
$P(k!|J_{\mu,k}(f)|>u)\le \frac{E(k!|J_{\mu,k}(f)|)^{2M}}{u^{2M}}$
with an appropriate choice of the parameter~$M$.
 
Proposition B states that in the case $M_0\le M\le\e n\sigma^2$
the inequality
$$
E\(n^{-k/2}k!I_{n,k}(f)\)^{2M}\le E((1+\beta(\e))\sigma\eta^k)^{2M}
$$
holds with a standard normal random variable $\eta$ and a function
$\beta(\e)$, $0\le\e\le1$ such that $\beta(\e)\to0$ if $\e\to0$, and
$\beta(\e)\le C$ with some universal constant $C=C(k)$ depending
only on the order~$k$ of the $U$-statistic for all $0\le\e\le1$.
This means that certain high but not too high moments of
$n^{-k/2}k!I_{n,k}(f)$ behave similarly to the moments of
$k!J_{\mu,k}(f)$. As a consequence, we can prove a similar, but
slightly weaker estimate for the distribution of
$n^{-k/2}k!I_{n,k}(f)$ as for the distribution of $k!J_{\mu,k}(f)$.
Theorem~3 contains the result we can get about the distribution of
$I_{n,k}(f)$ by means of these moment estimates.
 
The proof of Proposition~A is based on a corollary of a most
important result about Wiener--It\^o integrals called the diagram
formula. This result enables us to rewrite the product of
Wiener--It\^o integrals as a sum of Wiener--It\^o integrals of
different order.  It got the name `diagram formula' because
the kernel functions of the Wiener--It\^o integrals appearing in the
sum representation of the product of Wiener--It\^o integrals
are defined with the help of certain diagrams. As the expectation of
a Wiener--It\^o integral of order $k$ is zero for all $k\ge1$ the
expectation of the product equals the sum of the constant terms (i.e.
of the integrals of order zero) in the diagram formula. We shall see
that Proposition~A can be proved relatively simply by means of this
corollary of the diagram formula.
 
We shall also see that there is a version of the diagram formula
which enables us to express the product of degenerate $U$-statistics
as a sum of degenerate $U$-statistics of different order.
Proposition~B can be proved by means of this version of the diagram
formula similarly to the proof of Proposition~A. The main difference
between their proof is that in the case of the diagram formula for
degenerate $U$-statistics some new diagrams also appear, and their
contribution also has to be estimated. It will be shown
that if not too high moments of $U$-statistics are calculated by
means of this new version of the diagram formula, then the
contribution of the new diagrams is not too large.
 
The proof of formula (1.2) in [1] also contains the proof of the
inequality
$$
E\(n^{-k/2}k!I_{n,k}(f)\)^{2M}\le C^{2M}M^{kM} \tag2.4
$$
with some appropriate constant $C=C(k)$ for $M\le n\sigma^2$ in an
implicit way. This estimate is sufficient to the prove relation
(1.2), but insufficient to prove Theorem~3. In this case we
 need such a sharpened version of  inequality (2.4) which contains
an asymptotically optimal constant $C$ if $M\le \e n\sigma^2$ with
a small coefficient $\e>0$. But the method of paper~[1] is not
strong enough to prove such a sharpened version of~(2.4).
 
One reason for this weakness of the method of paper~[1] is that it
applies a consequence of Borell's inequality which does not give a
sharp inequality. Nevertheless, this inequality could be improved.
(See my paper~[9].) Another problem is that the proof in~[1]
contains a decoupling argument of paper~[2]. This argument which is
needed to apply a multivariate version of the Marcinkiewicz--Zygmund
inequality also weakens the universal constants in formula~(1.2).
This difficulty could also be overcome by some clever tricks.
But the application of the Marcinkiewicz--Zygmund inequality does
not allow to prove relation (2.4) with an optimal constant~$C$. The
proof of this inequality is based on a symmetrization argument which
implies in particular, that the moments of $n^{-k/2}k!I_{n,k}(f)$
are bounded by the (same) moments of a random variable with a
constant times greater variance. The influence of this too large
variance is inherited in all subsequent estimates, and as a
consequence, a method applying a symmetrization argument cannot
yield the estimate (2.4) with a sharp constant~$C$.
 
\beginsection 3. The proof of Theorem 1.
 
To formulate the corollary of the diagram formula we need in the
proof of Proposition~A first I introduce some notations.
 
Let us have a $\sigma$-finite measure $\mu$ together with a white
noise $\mu_W$ with counting measure $\mu$ on $(X,\Cal X)$. Let us
consider $L$ real valued functions $f_l(x_1,\dots,x_{k_l})$ on
$(X^{k_l},\Cal X^{k_l})$ such that
$\int f^2_l(x_1,\dots,x_{k_l})\mu(\,dx_1)\dots\mu(\,dx_{k_l})<\infty$,
$1\le l\le L$. Let us introduce the Wiener--It\^o integrals
$k_l!J_{\mu,k_l}(f_l)=\int f_l(x_1,\dots,x_{k_l})
\mu_W(\,dx_1)\dots\mu_W(\,dx_{k_l})$, $1\le l\le L$, and describe
how the expected value $E\(\prodd_{l=1}^L k_l!J_{\mu,k_l}(f_l)\)$
can be calculated by means of the diagram formula.
 
For this goal let us introduce the following notations. Put
$$
F(x_{(l,j)}, 1\le l\le L,\,1\le j\le k_l)=\prod_{l=1}^L
f_l(x_{(l,1)},\dots,x_{(l,k_l)}), \tag3.1
$$
and define a class of diagrams $\Gamma(k_1,\dots,k_L)$ in the
following way: Each diagram $\gamma\in\Gamma(k_1,\dots,k_L)$ is a
(complete, undirected) graph with vertices $(l,j)$, $1\le l\le L$,
$1\le j\le k_l$, and we shall call the set of vertices $(l,j)$ with
a fixed index $l$ the $l$-th row of a graph
$\gamma\in\Gamma(k_1,\dots,k_L)$. The graphs
$\gamma\in\Gamma(k_1,\dots,k_L)$ will have edges with the following
properties. Each edge connects vertices $(l,j)$ and $(l',j')$ from
different rows, i.e. $l\neq l'$ for the end-points of an edge. From
each vertex there starts exactly one edge. $\Gamma(k_1,\dots,k_L)$
contains all graphs $\gamma$ with such properties. If there is no
such graph, then $\Gamma(k_1,\dots,k_L)$ is empty.
 
Put $2N=\summ_{l=1}^Lk_l$. Then each $\gamma\in \Gamma(k_1,\dots,k_L)$
contains exactly $N$ edges. If an edge of the diagram $\gamma$
connects some vertex $(l,j)$ with some other vertex $(l',j')$,
$l'>l$, then we call $(l',j')$ the lower end-point of this edge,
and we denote the set of lower end-points of $\gamma$ by $\Cal
A_\gamma$ which has $N$ elements. Let us also introduce the
following function $\alpha_\gamma$ on the vertices of $\gamma$. Put
$\alpha_\gamma(l,j)=(l,j)$ if $(l,j)$ is the lower end-point of an
edge, and $\alpha_\gamma(l,j)=(l',j')$ if $(l,j)$ is connected with
the point $(l'j')$ by an edge of $\gamma$, and $(l',j')$ is the
lower end-point of this edge. Then we define the function
$$
\bar F_\gamma(x_{(l,j)}, \; (l,j)\in \Cal A_\gamma)=
F(x_{\alpha_\gamma(l,j)}, 1\le l\le L,\,1\le j\le k_l) \tag3.2
$$
with the function $F$ introduced in (3.1), i.e. we replace the
argument $x_{(l,j)}$ by $x_{(l',j')}$ in the function $F$ if $(l,j)$
and $(l',j')$ are connected by an edge in $\gamma$, and $l'>l$.
Then we enumerate the lower end-points somehow, and define the function
$B_\gamma(r)$, $1\le r\le N$, such that $B_\gamma(r)$ is the $r$-th
lower end-point of the diagram $\gamma$. Write
$$
F_\gamma(x_1,\dots,x_N)=\bar F_\gamma(x_{B_\gamma(r)},\; 1\le r\le N)
$$
and
$$
F_\gamma=\int\cdots\int
F_\gamma(x_1,\dots,x_N)\mu(\,dx_1)\dots\mu(\,dx_N)
\quad\text{for all }  \gamma\in \Gamma(k_1,\dots,k_L). \tag3.3
$$
(The function $F_\gamma(x_1,\dots,x_N)$  depends on the enumeration
of the lower end-points of the diagram $\gamma$, but its integral
$F_\gamma$ is independent of it.)
 
We shall need the following corollary of the diagram formula.
\medskip\noindent
{\bf Theorem A.} {\it With the above introduced notation
$$
E\(\prodd_{l=1}^L k_l!J_{\mu,k_l}(f_l)\)=\sum_{\gamma\in
\Gamma(k_1,\dots,k_L)} F_\gamma. \tag3.4
$$
(If $\Gamma(k_1,\dots,k_L)$ is empty, then the expected value of
the above product of random integrals equals zero.) Beside this
$$
F_\gamma^2\le \prod_{l=1}^L\int
 f^2_l(x_1,\dots,x_{k_l})\mu(\,dx_1)\dots\mu(dx_{k_l})
\quad\text{for all } \gamma\in\Gamma(k_1,\dots,k_L). \tag3.5
$$
}\medskip
The proof of Theorem A can be found in Corollary~5.4 of~[7] or~[6].
The result of~[7] actually deals with a different version of
Wiener-It\^o integrals where their `Fourier transforms' are
considered, and we integrate not with respect to a white noise, but
with respect to its `Fourier transform'. The results obtained for
such integrals are actually equivalent to the result formulated in
Theorem~A. I formulated Theorem~A in the present form because
generally this version of Wiener--It\^o integrals is applied in the
literature, and it can be compared better with the diagram formula for
the product of degenerate $U$-statistics applied in this paper.
Paper~[6] contains the diagram
formula for the version of Wiener--It\^o integrals considered in this
paper. The result of Theorem~A which is not contained explicitly
in~[6] can be deduced from the diagram formula proved in~[6] in the
same (simple) way as Corollary~5.4 is proved in~[7]. Now we turn to
the proof of Proposition~A.
\medskip\noindent
{\it Proof of Proposition A.}\/ Proposition A can be simply proved
with the help of Theorem~A if we apply it with $L=2M$, and the
functions $f_l(x_1,\dots,x_{k_l})=f(x_1,\dots,x_k)$ for all $1\le
l\le 2M$. Then Theorem~A yields that
$$
E\(k!J_{\mu,k}(f)^{2M}\)\le\( \int
f^2(x_1,\dots,x_k)\mu(\,dx_1)\dots\mu(dx_k)\)^M|\Gamma_{2M}(k)|,
$$
where $|\Gamma_{2M}(k)|$ denotes the number of diagrams $\gamma$
in $\Gamma(\underbrace{k,\dots,k}_{2M\text{ times}})$. Thus
to complete the proof of Proposition~A it is enough to show that
$|\Gamma_{2M}(k)|\le 1\cdot3\cdot5\cdots(2kM-1)$. But this can be
seen simply with the help of the following observation. Let
$\bar\Gamma_{2M}(k)$ denote the class of all graphs with vertices
$(l,j)$, $1\le l\le 2M$, $1\le j\le k$, such that from all vertices
$(l,j)$ exactly one edge starts, all edges connect different
vertices, but we also allow edges connecting vertices $(l,j)$ and
$(l,j')$ with the same first coordinate~$l$. Let
$|\bar\Gamma_{2M}(k)|$ denote the number of graphs in
$\bar\Gamma_{2M}(k)$. Then clearly
$|\Gamma_{2M}(k)|\le|\bar\Gamma_{2M}(k)|$. On the other hand,
$|\bar\Gamma_{2M}(k)|=1\cdot3\cdot5\cdots(2kM-1)$. Indeed, let us
list the vertices of the graphs from $\bar\Gamma_{2M}(k)$ in an
arbitrary way. Then the first vertex can be paired with another
vertex in $2kM-1$ way, after this the first vertex from which no
edge starts can be paired with $2kM-3$ vertices from which no edge
starts. By following this procedure the next edge can be chosen
$2kM-5$ ways, and by continuing this calculation we get the desired
formula.
\medskip\noindent
{\it Proof of Theorem 1}. By Proposition~A, formula (2.2) and the
Markov inequality we have
$$
P\(|k!J_{\mu,k}(f)|>u\)\le\frac{E\(k!J_{\mu,k}(f)\)^{2M}}{u^{2M}}
\le A\(\frac {2kM\sigma^{2/k}}{eu^{2/k}}\)^{kM} \tag3.6
$$
with some constant $A>\sqrt2$ if $M\ge M_0$ with some constant
$M_0=M_0(A)$, and $M$ is an integer.
 
Put $\bar M=\bar M(u)=\frac1{2k}\(\frac u\sigma\)^{2/k}$, and
$M=M(u)=[\bar M]$, where $[x]$ denotes the integer part of a real
number $x$. Choose some number $u_0$ such that $\frac1{2k}\(\frac
{u_0}\sigma\)^{2/k}\ge M_0+1$. Then we can apply relation (3.6) with
$M=M(u)$ for $u\ge u_0$, and it yields that
$$
\aligned
P\(|k!J_{\mu,k}(f)|>u\) &\le A\(\frac
{2kM\sigma^{2/k}}{eu^{2/k}}\)^{kM}\le e^{-kM} \le Ae^{k}e^{-k\bar M} \\
&=Ae^k\exp\left\{-\frac12
\(\frac u\sigma\)^{2/k}\right\} \quad\text{if } u\ge u_0.
\endaligned \tag3.7
$$
Relation (3.7) means that relation (1.6) holds for $u\ge u_0$ with
the pre-exponential coefficient $Ae^k$. By enlarging this
coefficient if it is needed we can guarantee that relation (1.6)
holds for all $u>0$. Theorem~1 is proved.
 
\beginsection 4. The diagram formula for the product of two degenerate
$U$-statistics.
 
To prove Proposition B we need a result analogous to Theorem~A
about the expectation of products of degenerate $U$-statistics.
To get such a result first we describe the product of two degenerate
$U$-statistics as the sum of degenerate $U$-statistics of different
order together with a good estimate on the $L_2$-norm of the kernel
functions in the sum representation. The proof of this result will
be given in the Appendix. We can get with the help of an inductive
procedure a generalization of this result. It yields a
representation of the product of several degenerate $U$-statistics
in the form of a sum of degenerate $U$-statistics which implies a
formula about the expected value of products of degenerate
$U$-statistics useful in the proof of Proposition~B. This
generalization will be discussed in the next section.
 
Let us have a sequence of iid. random variables $\xi_1,\xi_2,\dots$
with some distribution $\mu$ on a measurable space $(X,\Cal X)$
together with two functions $f(x_1,\dots,x_{k_1})$ and
$g(x_1,\dots,x_{k_2})$ on $(X^{k_1},\Cal X^{k_1})$ and
on $(X^{k_2},\Cal X^{k_2})$ respectively which are canonical
with respect to the probability measure $\mu$. We consider the
degenerate $U$-statistics $I_{n,k_1}(f)$ and $I_{n,k_2}(g)$
and express their normalized product
$k_1!k_2!n^{-(k_1+k_2)/2}I_{n,k_1}(f)I_{n,k_2}(g)$ as a sum of
(normalized) degenerate $U$-statistics. This product can be
written as a sum of $U$-statistics in a natural way, and then
by applying the Hoeffding decomposition for each of these
$U$-statistics as a sum of degenerate $U$-statistics we get the
desired representation of the product of two degenerate
$U$-statistics. The result we get in such a way will be presented
in Theorem~B. Before its formulation I introduce some notations.
 
To define the kernel functions of the $U$-statistics appearing in
the diagram formula for the product of two $U$-statistics first we
introduce
a class of objects $\Gamma(k_1,k_2)$ we shall call coloured diagrams.
We define graphs $\gamma\in\Gamma(k_1,k_2)$ that contain the vertices
$(1,1),(1,2),\dots,(1,k_1)$ which we shall call the first
row and $(2,1)\dots,(2,k_2)$ which we shall call the second row of
these graphs. From each vertex there starts zero or one edge, and
all edges connect vertices from different rows. All edges will get
a colour $+1$ or $-1$. $\Gamma(k_1,k_2)$  consists of all
$\gamma$ obtained in such a way which we shall call coloured
diagrams.
 
Given a coloured diagram $\gamma\in \Gamma(k_1,k_2)$ let
$B_u(\gamma)$ denote the set of upper end-points $(1,j)$ of the edges
of the graph $\gamma$, $B_{(b,1)}(\gamma)$ the set of lower end-points
$(2,j)$ of the edges of $\gamma$ with colour 1, and
$B_{(b,-1)}(\gamma)$ the set of lower end-points $(2,j)$ of the edges
of $\gamma$ with colour $-1$. (The letter `b' in the index was chosen
because of the word below.) Finally, let $Z(\gamma)$ denote the
set of edges with colour 1, $W(\gamma)$ the set of edges with
colour $-1$ of a coloured graph $\gamma\in\Gamma(k_1,k_2)$, and let
$|Z(\gamma)|$ and $|W(\gamma)|$ denote their cardinality.
 
Given two functions $f(x_1,\dots,x_{k_1})$ and $g(x_1,\dots,x_{k_2})$
let us define the function
$$
(f\circ g)(x_{(1,1)},\dots,x_{(1,k_1)},x_{(2,1)},\dots,x_{(2,k_2)})=
f(x_{(1,1)},\dots,x_{(1,k_1)})g(x_{(2,1)},\dots,x_{(2,k_2)}) \tag4.1
$$
 
Given a function $h(x_{u_1},\dots,x_{u_r})$ with coordinates in the
space $(X,\Cal X)$ (the indices $u_1,\dots,u_r$ are all different)
let us introduce  its transforms $P_{u_j}h$ and $O_{u_j}h$ by
the formulas
$$
P_{u_j}h(x_{u_l}\: u_l\in \{u_1,\dots,u_r\}\setminus \{u_j\})=
\int h(x_{u_1},\dots,x_{u_r})\mu(\,dx_{u_j}), \quad 1\le j\le r,
\tag4.2
$$
and
$$
Q_{u_j}h(x_{u_1},\dots,x_{u_r})=h(x_{u_1},\dots,x_{u_r})-
\int h(x_{u_1},\dots,x_{u_r})\mu(\,dx_{u_j}), \quad 1\le j\le r.
\tag4.3
$$
At this point I started to apply a notation which may seem to be
too complicated, but I think that it is more appropriate in the
further discussion. Namely, I started to apply a rather
general enumeration $u_1,\dots,u_r$ of the arguments of the
functions we are working with instead of their simpler enumeration
with indices $1,\dots,r$. But in the further discussion there will
appear an enumeration of the arguments by pairs of integers $(l,j)$
in a natural way, and I found it simpler to work with such an
enumeration than to reindex our variables all the time. Let me
remark in particular that this means that the definition
of the $U$-statistic with a kernel function $f(x_1,\dots,x_k)$
given in formula (1.1) will appear sometimes in the following more
complicated, but actually equivalent form: We shall work with kernel
function $f(x_{u_1},\dots,x_{u_k})$ instead of $f(x_1,\dots,x_k)$,
the random variables $\xi_j$ will be indexed by $u_s$, i.e. to the
coordinate $x_{u_s}$ we shall put the random variables
$\xi_{j_{u_s}}$ with indices $1\le j_{u_s}\le n$, and in the new
notation formula (1.1) will look like
$$
I_{n,k}(f)=\frac1{k!}\summ\Sb 1\le j_{u_s}\le n,\; s=1,\dots, k\\
j_{u_s}\neq j_{u'_{s}} \text{ if } u_s\neq u'_{s}\endSb
f\(\xi_{j_{u_1}},\dots,\xi_{j_{u_k}}\).   \tag$1.1'$
$$

Let us define for all coloured diagrams $\gamma\in\Gamma(k_1,k_2)$
the function $\alpha_\gamma(1,j)$, $1\le j\le k_1$, on the vertices
of the first row of $\gamma$ as $\alpha_\gamma(1,j)=(1,j)$ if no edge
starts from $(1,j)$, and $\alpha_\gamma(1,j)=(2,j')$ if an edge of
$\gamma$ connects the vertices $(1,j)$ and $(2,j')$. Given two
functions $f(x_1,\dots,x_{k_1})$ and $g(x_1,\dots,x_{k_2})$
together with a coloured diagram $\gamma\in\Gamma(k_1,k_2)$ let us
introduce, with the help of the above defined function
$\alpha_\gamma(\cdot)$ and $(f\circ g)$ introduced in (4.1) the
function
$$
\aligned
&\overline{(f\circ g)}_\gamma
(x_{(1,j)},x_{(2,j')},\,j\in\{1,\dots,k_1\}\setminus
B_u(\gamma), \, 1\le j'\le k_2) \\
&\qquad=(f\circ g)
(x_{\alpha_\gamma(1,1)},\dots,x_{\alpha_\gamma(1,k_1)},x_{(2,1)},
\dots,x_{(2,k_2)}).
\endaligned \tag4.4
$$
(In words, we take the function $(f\circ g)$, and if there is an
edge of $\gamma$ starting from a vertex $(1,j)$, and it connects
this vertex with the vertex $(2,j')$, then the argument $x_{(1,j)}$
is replaced by the argument $x_{(2,j')}$ in this function.) Let us
also introduce the function
$$
\aligned
&(f\circ g)_\gamma
\(x_{(1,j)},x_{(2,j')},\,j\in\{1,\dots,k_1\}\setminus
B_u(\gamma), \, j'\in \{1,\dots,k_2\}\setminus B_{(b,1)}\)\\
&\qquad=\prod_{(2,j')\in B_{(b,1)}(\gamma)} P_{(2,j')}
\prod_{(2,j')\in B_{(b,-1)}(\gamma)} Q_{(2,j')} \\
&\qquad\qquad\qquad \overline{(f\circ g)}_\gamma
\(x_{(j,1)},x_{(j',2)},\;j\in\{1,\dots,k_1\}\setminus
B_u(\gamma), \; 1\le j'\le k_2\).
\endaligned \tag4.5
$$
(In words, we take the function $\overline{(f\circ g)}_\gamma$ and
for such indices $(j',2)$ of the graph $\gamma$ from which an edge
with colour 1 starts we apply the operator $P_{(2,j')}$ introduced
in formula (4.2) and for those indices $(2,j')$ from which an edge
with colour $-1$ starts we apply the operator $Q_{(2,j')}$ defined
in formula~(4.3).) Let us also remark that the operators
$P_{(2,j')}$ and $Q_{(2,j')}$ are exchangeable for different
indices $j'$, hence it is not important in which order we apply the
operators $P_{(2,j')}$ and $Q_{(2,j')}$ in formula (4.5).
 
In the definition of the function $(f\circ g)_\gamma$ those arguments
$x_{(2,j')}$ of the function $\overline{(f\circ g)}_\gamma$ which
are indexed by a pair $(2,j')$ from which an edge of colour~1 of the
coloured diagram $\gamma$ starts will disappear, while the arguments
indexed by a pair $(2,j')$ from which an edge of colour $-1$ of the
coloured diagram $\gamma$ starts will be preserved.
Hence the number of arguments in the function $(f\circ g)_\gamma$
equals $k_1+k_2-2|B_{(b,1)}(\gamma)|-|B_{(b,-1)}(\gamma)|$, where
$|B_{(b,1)}(\gamma)|$ and $|B_{(b,-1)}(\gamma)|$ denote the
cardinality of the lower end-points of the edges of the coloured
diagram $\gamma$ with colour 1 and $-1$ respectively, In an
equivalent form we can say that the number of arguments of
$(f\circ g)_\gamma$ equals $k_1+k_2-(2|Z(\gamma)|+|W(\gamma)|)$.
 
Now we are in the position to formulate the diagram formula for
the product of two degenerate $U$-statistics.
\medskip\noindent
{\bf Theorem B.} {\it Let us have a sequence of iid. random
variables $\xi_1,\xi_2,\dots$ with some distribution $\mu$ on some
measurable space $(X,\Cal X)$ together with two bounded, canonical
functions $f(x_1,\dots,x_{k_1})$ and $g(x_1,\dots,x_{k_2})$ with
respect to the probability measure~$\mu$ on the spaces $(X^{k_1},
\Cal X^{k_1})$ and $(X^{k^2},\Cal X^{k_2})$. Let us introduce the
class of coloured diagrams $\Gamma(k_1,k_2)$ defined above together
with the functions $(f\circ g)_\gamma$ defined in formulas
(4.1)---(4.5).
 
For all $\gamma\in\Gamma$ the function $(f\circ g)_\gamma$ is
canonical with respect to the measure $\mu$ with
$k(\gamma)=k_1+k_2-(2|Z(\gamma)|+|W(\gamma)|)$
arguments, where $|Z(\gamma)|$ denotes the number of edges with
colour 1 and $|W(\gamma)|$ the number of edges with colour $-1$ of
the coloured diagram~$\gamma$. The product of the degenerate
$U$-statistics $I_{n,k_1}(f)$ and $I_{n,k_2}(g)$, $n\ge
\max(k_1,k_2)$, defined in (1.1) satisfies the identity
$$ \allowdisplaybreaks
\align
&k_1!k_2!n^{-(k_1+k_2)/2}I_{n,k_1}(f)I_{n,k_2}(g)\\
&\qquad=\sumn_{\gamma\in\Gamma(k_1,k_2)}
\frac{\prodd_{j=1}^{|Z(\gamma)|}
(n-(k_1+k_2)+|W(\gamma)|+|Z(\gamma)|+j)}{n^{|Z(\gamma)|}} \tag4.6 \\
&\qquad\qquad\qquad \; n^{-|W(\gamma)|/2}\cdot
k(\gamma)!n^{-k(\gamma)/2}I_{n,k(\gamma)}((f\circ g)_\gamma),
\endalign
$$
where $\sum^{\prime(n)}$ means that summation is taken only for
such coloured diagrams $\gamma\in\Gamma(k_1,k_2)$ which satisfy
the inequality $k_1+k_2-(|Z(\gamma)|+|W(\gamma)|)\le n$, and
$\prodd_{j=1}^{|Z(\gamma)|}$ equals 1 in the case $|Z(\gamma)|=0$.
 
The $L_2$-norm of the functions $(f\circ g)_\gamma$ is defined by
the formula
$$ \allowdisplaybreaks
\align
\|(f\circ g)_\gamma\|_2^2&=
\int (f\circ g)_\gamma^2
(x_{(1,j)},x_{(2,j')},\,j\in\{1,\dots,k_1\}\setminus
B_u(\gamma), \, j'\in \{1,\dots,k_2\}\setminus B_{(b,1)}) \\
&\qquad \prod_{(1,j)\: j\in\{1,\dots,k_1\}\setminus
B_u(\gamma)}\mu(\,dx_{(1,j)})
\prod_{(2,j')\: j'\in\{1,\dots,k_2\}\setminus B_{(b,1)}}
\mu(\,dx_{(2,j')}).
\endalign
$$
If $W(\gamma)=0$, then the inequality
$$
\|(f\circ g)_\gamma\|_2\le \|f\|_2\|g\|_2
\quad \text{if }\; |W(\gamma)|=0. \tag4.7
$$
holds. In the general case we can say that if the functions $f$
and $g$ satisfy formula (1.3), then also the inequality
$$
\|(f\circ g)_\gamma\|_2\le 2^{|W(\gamma)|}\min(\|f\|_2,\|g\|_2)
\tag4.8
$$
holds. Relations (4.7) and (4.8) remain valid even if we drop the
condition that the functions $f$ and $g$ are canonical.}\medskip
 
Relations (4.7) and (4.8) mean in particular, that we have a better
estimate for $\|(f\circ g)_\gamma\|_2$ in the case when the coloured
diagram $\gamma$ contains no edge with colour $-1$, i.e.
$|W(\gamma)|=0$, than in the case when it contains at least one
edge with colour $-1$.
 
Let us understand how we define those terms at the right-hand side
of (4.6) for which $k(\gamma)=0$. In this case $(f\circ g)_\gamma$
is a constant, and to make formula (4.6) meaningful we have to
define the term $I_{n,k(\gamma)}((f\circ g)_\gamma)$ also in this
case. The following convention will be used. A constant $c$ will be
called a degenerate $U$-statistic of order zero, and we define
$I_{n,0}(c)=c$.
 
Theorems~B can be considered as a version of the result of paper~[8],
where a similar diagram formula was proved about multiple random
integrals with respect to normalized empirical measures. Degenerate
$U$-statistics can also be presented as such integrals with special,
canonical kernel functions. Hence there is a close relation between
the results of this paper and~[8]. But there are also some essential
differences. For one part, the diagram formula for multiple random
integrals with respect to normalized empirical measures is simpler
than the analogous result about the product of degenerate
$U$-statistics, because the kernel functions in these integrals need
not be special, canonical functions. On the other hand, the
 diagram formula for degenerate $U$-statistics yields a simpler
formula about the expected value of the product of degenerate
$U$-statistics, because the expected value of a degenerate
$U$-statistic equals zero, while the analogous result about multiple
random integrals with respect to normalized empirical measures does
not hold. Another difference between this paper and~[8] is that
here I worked out a new notation which, I hope, is more transparent.
 
\beginsection 5. The diagram formula for the product of several
degenerate $U$-statistics.
 
We can also express the product of more than two degenerate
$U$-statistics in the form of sums of degenerate $U$-statistics
by applying Theorem~B recursively. We shall present this result in
Theorem~B$'$ and prove it together with an estimate about the
$L_2$-norm of the kernel functions of the degenerate~$U$-statistics
appearing in Theorem~B$'$. This estimate will be given in Theorem~C.
Since the expected value of all degenerate $U$-statistics of order
$k\ge1$ equals zero, the representation of the product of
$U$-statistics in the form of a sum of degenerate $U$-statistics
implies that the expected value of this product equals the sum of
the constant terms in this representation. In such a way we get a
version of Theorem~A for the expected value of a product of
degenerate $U$-statistics which together with Theorem~C will be
sufficient to prove Proposition~B. But the formula we get in this
way is more complicated than the analogous diagram formula for
products of Wiener--It\^o integrals. To overcome this difficulty we
have to work out a good ``book-keeping method''.
 
Let us have a sequence of iid. random variables $\xi_1,\xi_2,\dots$
taking values on a measurable space $(X,\Cal X)$ with some
distribution~$\mu$, and consider $L$ functions $f_l(x_1,\dots,x_{k_l})$
on the measure spaces $(X^{k_l},\Cal X^{k_l})$, $1\le l\le L$,
canonical with respect to the measure~$\mu$. We want to represent
the product of $L\ge2$ normalized degenerate $U$-statistics
$n^{-k_l/2}k_l!I_{n, k_l}(f_{k_l})$ in the form of a sum of
degenerate $U$-statistics similarly to Theorem~B. For this goal I
define a class of coloured diagrams $\Gamma(k_1,\dots,k_L)$ together
with some canonical functions
$F_\gamma=F_\gamma(f_{k_1},\dots,f_{k_L})$ depending on the
diagrams $\gamma\in\Gamma(k_1,\dots,k_L)$ and the functions
$f_l(x_1,\dots,x_{k_l})$, $1\le l\le L$.
 
The coloured diagrams will be graphs with vertices $(l,j)$ and
$(l,j,C)$, $1\le l\le L$, $1\le j\le k_l$, and edges between some of
these vertices which will get either colour 1 or colour $-1$. The
set of vertices $\{(l,j),(l,j,C),\,1\le j\le k_l\}$ will be called
the $l$-th row of the diagrams. (The vertices $(l,j,C)$ are
introduced, because it turned out to be useful to take a copy
$(l,j,C)$ of some vertices $(l,j)$. The letter $C$ was just chosen
to indicate that it is a copy.) From all vertices there starts either
zero or one edge, and edges may connect only vertices in different
rows. We shall call all vertices of the form $(l,j)$ permissible,
and beside this some of the vertices $(l,j,C)$ will also be called
permissible. Those vertices will be called permissible from which
some edge may start.
 
We shall say that an edge connecting two vertices $(l_1,j_1)$ with
$(l_2,j_2)$ or (a permissible) vertex $(l_1,j_1,C)$ with another
vertex $(l_2,j_2)$ such that $l_2>l_1$ is of level $l_2$, and
$(l_2,j)$ will be called the lower end-point of such an edge. (The
coloured diagrams we shall define contain only edges with lower
end-points of the form $(l,j)$.) We shall call the restriction
$\gamma(l)$ of the diagram $\gamma$ to level $l$ that part of a
diagram $\gamma$ which contains all of its vertices together with
those edges (together with their colours) whose levels are less
than or equal to $l$, and tells which of the vertices $(l',j,C)$
are permissible for $1\le l'\le l$. We shall define the diagrams
$\gamma\in\Gamma(k_1,\dots,k_L)$ inductively by defining their
restrictions $\gamma(l)$ to level $l$ for all $l=1,2,\dots,L$.
Those diagrams $\gamma$ will belong to $\Gamma(k_1,\dots,k_L)$
whose restrictions $\gamma(l)$ can be defined through the
following procedure for all $l=1,2,\dots,L$.
 
The restriction $\gamma(1)$ of a diagram $\gamma$ to level 1
contains no edges, and no vertex of the form $(1,j,C)$, $1\le j\le
k_1$, is permissible. If we have defined the restrictions
$\gamma(l-1)$ for some $2\le l\le L$, then those diagrams will be
called restrictions $\gamma(l)$ at level~$l$ which can be obtained
from a restriction $\gamma(l-1)$ in the following way: Take the
vertices $(l,j)$, $1\le j\le k_l$, from the $l$-th row and from each
of them either no edge starts or one edge starts which gets either
colour~1 or colour~$-1$. The other end-point must be such a vertex
$(l',j')$ or a permissible vertex $(l',j',C)$ with some $1<l'<l$
which is not an end-point of a vertex in $\gamma(l-1)$, and
naturally such a vertex can be connected only with one of the
vertices $(l,j)$, $1\le j\le k_l$. We define $\gamma(l)$ first by
adjusting the coloured edges constructed in the above way to the
(coloured) edges of $\gamma(l-1)$, and the set of permissible
vertices in $\gamma(l)$ will contain beside the permissible
vertices of $\gamma(l-1)$ and the vertices $(l,j)$, $1\le j\le k_l$,
those vertices $(l,j,C)$ for which $(l,j)$ is the lower end-point of
an edge with colour $-1$ in $\gamma(l)$. $\Gamma(k_1,\dots,k_L)$ will
consist of all coloured diagrams $\gamma=\gamma(L)$ obtained in such
a way.
 
Given a coloured diagram $\gamma\in\Gamma(k_1,\dots,k_L)$ we shall
define recursively some (canonical) functions $F_{l,\gamma}$ with the
help of the functions $f_1,\dots,f_l$ together with some constants
$J_n(l,\gamma)$ for all $1\le l\le L$ in the way suggested by
Theorem~B. Then we put $F_\gamma=F_{L,\gamma}$ and give the desired
representation of the product of the degenerate $U$-statistics
with the help of $U$-statistics with kernel functions $F_\gamma$
and constants $J_n(l,\gamma)$, $\gamma\in\Gamma(k_1,\dots,k_L)$,
$1\le l\le L$.
 
Let us fix some coloured diagram $\gamma\in\Gamma(k_1,\dots,k_L)$
and introduce the following notations: Let $B_{(b,-1)}(l,\gamma)$
denote the set of lower end-points of the form $(l,j)$ of edges with
colour $-1$ and $B_{(b,1)}(l,\gamma)$ the set of lower end-points of
the form $(l,j)$ with colour~$1$. Let $U(l,\gamma)$ denote the
set of those permissible vertices $(l',j)$ and $(l',j,C)$ with
$l'\le l$ from which no edge starts in the restriction $\gamma(l)$ of
the diagram $\gamma$ to level $l$, i.e. either no edge starts from
this vertex, or if some edge starts from it, then its other end-point
is a vertex $(l',j)$ with $l'>l$. Beside this, given some integer
$1\le l_1<l$ let $U(l,l_1,\gamma)$ denote the restriction of
$U(l,\gamma)$ to its first $l_1$ rows, i.e. $U(l,l_1,\gamma)$
consists of those vertices $(l',j)$ and $(l,j',C)$ which are
contained in $U(l,\gamma)$, and $l'\le l_1$. We shall define the
functions $F_l(\gamma)$ with arguments of the form $x_{(l',j)}$ and
$x_{(l',j,C)}$ with $(l',j)\in U(l,\gamma)$ and $(l',j,C)\in
U(l,\gamma)$ together with some constants $J_n(l,\gamma)$. For this
end put first
$$
F_{1,\gamma}(x_{(1,1)},\dots,x_{(k_1,1)})
=f_1(x_{(1,1)},\dots,x_{(k_1,1)}). \tag5.1
$$
To define the function $F_{l,\gamma}$ for $l\ge2$ first we introduce
a function $\alpha_{l,\gamma}(\cdot)$ on the set of vertices in
$U(l-1,\gamma)$ in the following way. If a vertex $(l',j')$ or
$(l',j',C)$ in $U(\gamma,l-1)$ is such that it is connected to no
vertex $(l,j)$, $1\le j\le k_l$, then
$\alpha_{l,\gamma}(l',j')=(l',j')$,
$\alpha_{l,\gamma}(l',j',C)=(l',j',C)$ and if $(l',j')$ is connected
to a vertex $(l,j)$, then $\alpha_{l,\gamma}(l',j')=(l,j)$,
if $(l',j',C)$ is connected with a vertex $(l,j)$, then
$\alpha_{l,\gamma}(l',j',C)=(l,j)$. We define, similarly to the
formula (4.4) the functions
$$
\aligned
&\bar F_{l,\gamma}(x_{(l',j')},x_{(l',j',C)},\;(l',j')\text{ and }
(l',j',C)\in U(l,l-1,\gamma),\; x_{(l,j)},\, 1\le j\le k_l)\\
&\qquad =F_{l-1,\gamma}
(x_{\alpha_{l,\gamma}(l',j')},x_{\alpha_{l,\gamma}(l',j',C)},\,
(l',j')\text{ and } (l',j',C)\in U(l-1,\gamma))\\
&\qquad\qquad\qquad\qquad f_l(x_{(l,1)},\dots,x_{(l,k_l)}),
\endaligned \tag5.2
$$
i.e. we take the function $F_{l-1,\gamma}\circ f_l$ and replace the
arguments of this function indexed by such a vertex of $\gamma$
which is connected by an edge with a vertex in the $l$-th row of
$\gamma$ by the argument indexed with the lower end-point of this
edge.
 
Then we define with the help of the operators $P_{u_j}$ and $Q_{u_j}$
introduced in (4.2) and (4.3) the functions
$$
\aligned
&\bar{\bar F}_{l,\gamma}(x_{(l',j')},x_{(l',j',C)},\;(l',j')
\text{ and }(l',j',C)\in U(l,l-1,\gamma),\\
&\qquad\qquad\qquad x_{(l,j)},\; j\in\{1,\dots,k_l\}\setminus
B_{(l,1)}(l,\gamma))\\
&\qquad=\prod_{(l,j)\in B_{(b,1)}(l,\gamma)} P_{(l,j)}
\prod_{(l,j)\in B_{(b,-1)}(l,\gamma)}  Q_{(l,j)} \\
&\qquad\qquad \bar F_{l,\gamma}(x_{(l',j')},x_{(l',j',C)},\;(l',j')
\text{ and } (l',j',C)\in U(l,l-1,\gamma),\;
x_{(l,j)},\; 1\le j\le k_l),
\endaligned \tag5.3
$$
similarly to the formula (4.5), i.e. we apply for the function $\bar
F_l(\gamma)$ the operators $P_{(l,j)}$ for those indices $(l,j)$
which are the lower end-points of an edge with colour 1 and the
operators $Q_{(l,j)}$ for those indices $(l,j)$ which are the lower
end-points of an edge with colour $-1$.
 
Finally we define the function $F_{l,\gamma}$ simply by reindexing
some arguments of the function $\bar{\bar F}_{l,\gamma}$ to get a
function which is indexed by the vertices in $U(l,\gamma)$. To this
end we define the function $A_{l,\gamma}(\cdot)$ on the set of
vertices $\{(l,j)\:(l,j)\in\{(l,1),\dots,(l,k_l)\}\setminus
B_{(b,1)}(l,\gamma)$
as $A_{l,\gamma}(l,j)=(l,j,C)$ if $(l,j)\in B_{(b,-1)}(l,\gamma)$,
and $A_{l,\gamma}(l,j)=(l,j)$ if
$(l,j)\in\{(l,1),\dots,(l,k_l)\}\setminus
(B_{(b,1)}(l,\gamma)\cup B_{(b,-1)}(l,\gamma))$. Then we put
$$
\aligned
&F_{l,\gamma}(x_{(l',j')},x_{(l',j',C)},\; (l',j')\text{ and }
(l',j',C)\in U(l,\gamma))\\
&\qquad=\bar{\bar F}_{l,\gamma}(x_{(l',j')}, x_{(l',j',C)},\;(l',j')
\text{ and } (l',j,C)\in U(l,l-1,\gamma),\\
&\qquad\qquad\qquad
x_{A_{l,\gamma}(l,j)},\; (l,j)\in\{(l,1),\dots,(l,k_l)\}\setminus
B_{(b,1)}(l,\gamma)).
\endaligned \tag5.4
$$
 
We define beside the functions $F_\gamma=F_{L,\gamma}$ the following
constants $J_n(l,\gamma)$, $1\le l\le L$: $J_n(1,\gamma)=1$, and
$$
J_n(l,\gamma)= \frac{\prodd_{j=1}^{|B_{(b,1)}(l,\gamma)|}
(n-(k_1+k_2)+|B_{(b,-1)}(l,\gamma)|+|B_{(b,1)}(l,\gamma)|+j)}
{n^{|Z(\gamma)|}},
\quad 2\le l\le L, \tag5.5
$$
if $|B_{(b,1)}(l,\gamma)|\ge1$, and $J_n(l,\gamma)=1$ if
$|B_{(b,1)}(l,\gamma)|=0$, where $|B_{(b,1)}(l,\gamma)|$ and
$|B_{(b,-1)}(l,\gamma)|$ denote the number of those edges in
$\gamma$ with colour 1 and with colour $-1$ respectively whose
lower end-point is in the $l$-th row of $\gamma$.
 
Now we can formulate the following generalization of Theorem~B.
\medskip\noindent
{\bf Theorem $\bold B'$.} {\it Let us have a sequence of iid. random
variables $\xi_1,\xi_2,\dots$ with some distribution $\mu$ on a
measurable space $(X,\Cal X)$ together with $L\ge2$ bounded functions
$f_l(x_1,\dots,x_{k_l})$ on the spaces $(X^{k_l},\Cal X^{k_l})$,
$1\le l\le L$, canonical with respect to the probability
measure~$\mu$.  Let us introduce the class of coloured diagrams
$\Gamma(k_1,\dots,k_L)$ defined above together with the functions
$F_\gamma=F_{L,\gamma}(f_1,\dots,f_L)$ defined in formulas
(5.1)---(5.4) and the constants $J_n(l,\gamma)$, $1\le l\le L$ given
in formula (5.5).
 
Put $k(\gamma(l))=\summ_{p=1}^l k_p -\summ_{p=2}^l
(2|B_{(b,1)}(p,\gamma)|+|B_{(b,-1)}(p,\gamma)|)$, where
$|B_{(b,1)}(p,\gamma)|$ denotes the number of lower end-points in
the $p$-th row of $\gamma$ with colour 1 and $|B_{(b,-1)}(p,\gamma)|$
is the number of lower end-points in the $p$-th row of $\gamma$ with
colour $-1$, $1\le l\le L$, and define $k(\gamma)=k(\gamma(L))$.
Then $k(\gamma(l))$ is the number of variables of the function
$F_{l,\gamma}$, $1\le l\le L$.
 
The functions $F_\gamma$ are canonical with respect to the measure
$\mu$ with $k(\gamma)$ variables, and the product of the degenerate
$U$-statistics $I_{n,k_l}(f)$, $n\ge \max\limits_{1\le l\le L} k_l$,
defined in (1.1) satisfies the identity
$$
\prod_{l=1}^L k_l! n^{-k_l/2}I_{n,k_l}(f_{k_l})=
\sumnl_{\gamma\in\Gamma(k_1,\dots,k_L)} \(\prod_{l=1}^L J_n(l,\gamma)\)
n^{-|W(\gamma)|/2}\cdot k(\gamma)!n^{-k(\gamma)/2}I_{n,k(\gamma)}
(F_\gamma), \tag5.6
$$
where $|W(\gamma)|=\summ_{l=2}^L |B_{(b,-1)}(l,\gamma)|$ is the
number of edges with colour $-1$ in the coloured diagram $\gamma$,
and $\sum^{\prime(n,\,L)}$ means that summation is taken
for those $\gamma\in\Gamma(k_1,\dots,k_L)$ which satisfy the relation
$k(\gamma(l-1))+k_{l}-(|B_{(b,1)}(l,\gamma)|+|B_{(b,-1)}(l,\gamma)|)
\le n$ for all $2\le l\le L$.
 
Let $\bar\Gamma(k_1,\dots,k_L)$ denote the class of those coloured
diagrams of $\Gamma(k_1,\dots,k_L)$ for which every permissible
vertex is the end-point of some vertex. A coloured diagram
$\gamma\in\Gamma(k_1,\dots,k_L)$ satisfies the relation
$\gamma\in\bar\Gamma(k_1,\dots,k_L)$ if and only if $k(\gamma)=0$.
In this case $F_\gamma$ is constant, and
$I_{n,k(\gamma)}(F_\gamma)=F_\gamma$.  For all other coloured
diagrams $\gamma\in \Gamma(k_1,\dots,k_L)$ $k(\gamma)\ge0$. The
identity
$$
E\(\prod_{l=1}^L k_l! n^{-k_l/2}I_{n,k_l}(f_{k_l})\)
= \sumnl_{\gamma\in\bar\Gamma(k_1,\dots,k_L)} \(\prod_{l=1}^L
J_n(l,\gamma)\) n^{-|W(\gamma)|/2}\cdot F_\gamma
\tag5.7
$$
holds.
} \medskip
 
Theorem~B$'$ can be deduced relatively simply from Theorem~B by
induction with respect to the number $L$ of the functions. Theorem~B
contains the results of Theorem~B$'$ in the case $L=2$. A simple
induction argument together with the formulas describing the
functions $F_{l,\gamma}$ by means of the functions $F_{l-1,\gamma}$
and $f_l$ and Theorem~B imply that all functions $F_\gamma$ in
Theorem~B$'$ are canonical. Finally, an inductive procedure with
respect to the number~$L$ of the functions $f_l$ shows that
relation~(5.6) holds. Indeed, by exploiting that formula (5.6) holds
for the product of the first $L-1$ degenerate $U$-statistics, then
multiplying this identity with the last $U$-statistic and applying
for each term at the right-hand side Theorem~B we get that
relation~(5.6) also holds for the product $L$ degenerate
$U$-statistics.
 
A simple inductive procedure with respect to $l$ shows that for
all $2\le l\le L$  the diagram $\gamma(l)$ contains
$k(\gamma(l))=\summ_{p=1}^l k_l-\summ_{p=2}^l
(2|B_{(b,1)}(p,\gamma)|+|B_{(b,-1)}(p,\gamma)|)$ permissible vertices
in its first $l$ rows which is not the end-point of an edge in
$\gamma(l)$. Since $\gamma$ has $\summ_{p=1}^L k_l+\summ_{p=2}^L
|B_{(b,-1)}(p,\gamma)|)$ permissible vertices this identity
with $l=L$ implies that $k(\gamma)=0$ if and only if
$\gamma\in\bar\Gamma(k_1,\dots,k_L)$ with the class of coloured
diagrams $\bar\Gamma(k_1,\dots,k_L)$ introduced at the end of
Theorem~B$'$. Since $EI_{n,k}(f)=0$ for all degenerate
$U$-statistics of order $k\ge1$, the above property
and relation (5.6) imply identity (5.7).
 
In the proof of Proposition~B we shall also need an estimate
formulated in Theorem~C. It is a simple consequence of
inequalities~(4.7) and~(4.8) in Theorem~B.
\medskip\noindent
{\bf Theorem C.} {\it Let us have $L$ functions
$f_l(x_1,\dots,x_{k_l})$ on the spaces $(X^{k_l},\Cal X^{k_l})$,
$1\le l\le L$, which satisfy formulas (1.3) and (1.4) (if we
replace the index $k$ by index $k_l$ in these formulas), but these
functions need not be canonical. Let us take a coloured diagram
$\gamma\in\Gamma(k_1,\dots,k_L)$ and consider the function
$F_\gamma=F_{L,\gamma}(f_1,\dots,f_L)$ defined by formulas
(5.1)---(5.5). The $L_2$-norm of the function $F_\gamma$ (with
respect to a power of the measure $\mu$ to the space, where
$F_\gamma$ is defined) satisfies the inequality $\|F_\gamma\|_2\le
2^{|W(\gamma)|}\sigma^{(L-U(\gamma))}$, where $|W(\gamma)|$ denotes
the number of edges of  colour $-1$, and $U(\gamma)$ the number of
 rows which contain a lower vertex of colour $-1$ in the coloured
diagram~$\gamma$.}
\medskip\noindent
{\it Proof of Theorem C.} We shall prove the inequality
$$
\|F_{l,\gamma}\|_2\le 2^{|W(l,\gamma)|}
\sigma^{(l-U(l,\gamma))}\quad \text{for all } 1\le l\le L, \tag5.8
$$
where $|W(l,\gamma)|$ denotes the number of edges with colour 1,
and $U(l,\gamma)$ is the number of rows containing a lower point of
an edge with colour \ $-1$ in the coloured diagram $\gamma(l)$.
Formula (5.8) will be proved by means of induction with respect
to~$l$. It implies Theorem~C with the choice~$l=L$.
 
Relation (5.8) clearly holds for $l=1$. To prove this relation by
induction with respect to $l$ for all $1\le l\le L$ let us first
observe that $\sup 2^{-|W(l,\gamma)|}|F_{l,\gamma}|\le1$ for all
$1\le l\le L$. This relation can be simply checked by induction
with respect to~$l$.
 
If we know relation (5.8) for $l-1$, then it follows for $l$ from
relation~(4.7) if $|B_{(b,-1)}(l,\gamma)|=0$, that is if
there is no edge of colour \ $-1$ with lower end-point in the $l$-th
row. Indeed, in this case $\|F_{l,\gamma}(f_1,\dots,f_l)\|_2\le
\|F_{l-1,\gamma}\|_2 \|f_l\|_2\le
\|F_{l-1,\gamma}(f_1,\dots,f_{l-1})\|_2\cdot\sigma$,
$|W(l,\gamma)|=|W(l-1,\gamma)|$, and $U(l,\gamma)=U(l-1,\gamma)$.
Hence relation (5.8) holds in this case.
 
If $|B_{(b,-1)}(l,\gamma)|\ge1$, then we can apply formula (4.8)
for the expression $\|F_{l,\gamma}\|_2=\|\bar{\bar F}_{l,\gamma}\|_2=
\|(F_{l-1,\gamma}\circ f_l)_{\tilde\gamma(l)}\|_2$, where
$\tilde\gamma(l)$ is that coloured diagram with two rows whose
first row consists of the indices of the variables of the function
$F_{l-1,\gamma}$, its second row consists of the vertices $(l,j)$,
$1\le j\le k_l$, and $\tilde\gamma(l)$ contains the edges of
$\gamma$ between these vertices together with their colour. Then
relation (4.8) implies
that $\|F_{l,\gamma}\|_2\le2^{|B_{(b,-1)}|}\|F_{l-1,\gamma}\|_2
\le 2^{(|W(l-1,\gamma)|+|B_{(b,-1)}(l,\gamma)|)}
\sigma^{(l-1-U(l-1,\gamma))}$ if $|B_{(b,-1)}(l,\gamma)|\ge1$.
Beside this, $|W(l-1,\gamma)|+|B_{(b,-1)}(l,\gamma)|=|W(l,\gamma)|$,
and $l-1-U(l-1,\gamma)=l-U(l,\gamma)$ in this case. Hence
relation~(5.8) holds in this case, too.
 
\beginsection 6. The proof of Theorem 3.
 
First we prove Proposition B.
\medskip\noindent
{\it Proof of Proposition B.}\/ We shall prove relation (2.3) by
means of Theorem~C and identity (5.7) with the choice $L=2M$ and
$f_l(x_1,\dots,x_{k_l})=f(x_1,\dots,x_k)$ for all $1\le l\le 2M$.
We shall partition the class of coloured diagrams
$\gamma\in\Gamma(k,M)=\bar\Gamma(\underbrace{k,\dots,k}_{2M\text{
times}})$ with the property that all permissible vertices are the
end-points of some edge to classes $\Gamma(k,M,p)$, $1\le p\le M$,
in the following way: $\gamma\in\Gamma(k,M,p)$ for a coloured
diagram $\gamma\in\Gamma(M,k)$ if and only if it has $2p$
permissible vertices of the form $(l,j,C)$.
(A coloured diagram $\gamma\in \Gamma(k,M)$ has even number of
such vertices.)  First we prove the following estimate:
\medskip
{\narrower \noindent
There exists some constant $A=A(k)>0$ and threshold index
$M_0=M_0(k)$ such that for all $M\ge M_0$ and $0\le p\le kM$ the
cardinality $|\Gamma(k,M,p)|$ of the set $\Gamma(k,M,p)$ can be
bounded from above by
$A2^{2p}\binom{2kM}{2p}\(\frac2e\)^{kM}(kM)^{kM+p}$. \par}
 
\medskip
We can bound the number of coloured diagrams in $\Gamma(k,M,p)$ by
calculating first the number of choices of the $2p$ permissible
vertices from the $2kM$ vertices of the form $(l,j,C)$ which we
adjust to the $2kM$ permissible vertices $(l,j)$ and then by
calculating the number of such graphs whose vertices are the above
permissible vertices, and from all vertices there starts exactly one
edge. (Here we allow to connect vertices from the same row. Observe
that by defining the set of permissible vertices $(l,j,C)$ in a
coloured diagram $\gamma$ we also determine the colouring of its
edges.) Thus we get by using the argument at the beginning of
Proposition~A that $|\Gamma(k,M,p)|$ can be bounded from above by
$\binom{2kM}{2p}1\cdot3\cdot5\cdots(2kM+2p-1)=\binom{2kM}{2p}
\frac{(2kM+2p)!}{2^{kM+p}(kM+p)!}$. We can write by the Stirling
formula, similarly to relation (2.2) that
$\frac{(2kM+2p)!}{2^{kM+p}(kM+p)!}\le A\(\frac2e\)^{kM+p}(kM+p)^{kM+p}$
with some constant $A>\sqrt2$ if $M\ge M_0$ with some $M_0=M_0(A)$.
Since $p\le kM$ we can write $(kM+p)^{kM+p}\le (kM)^{kM}\(1+\frac
p{kM}\)^{kM}(2kM)^p\le (kM)^{kM+p}e^p2^p$. The above inequalities
imply that
$$
|\Gamma(k,M,p)|\le A\binom{2kM}{2p}\(\frac2e\)^{kM}(kM)^{kM+p}2^{2p}
\quad\text{if }M\ge M_0, \tag6.1
$$
as we have claimed.
 
Observe that for $\gamma\in\Gamma(k,M,p)$ the quantities introduced
in the formulation of Theorems~B$'$ and~$C$ satisfy the relations
$|W(\gamma)|=2p$, $|F_\gamma|=\|F_\gamma\|_2$ and $U(\gamma)\le
|W(\gamma)|=2p$. Hence by Theorem~C we have
$n^{-|W(\gamma)|/2}|F_\gamma| \le 2^p n^{-p}\sigma^{2M-U(\gamma)}\le
2^p \(n\sigma^2\)^{-p}\sigma^{2M}\le\eta^{p}2^p(kM)^{-p}
\sigma^{2M}$ if $kM\le \eta n\sigma^2$ and $\sigma^2\le1$.
 
This estimate together with relation (5.7) and the fact that
the constants $J_n(l,\gamma)$ defined in (5.5) are bounded by 1
imply that for $kM\le \eta n\sigma^2$
$$
E\(n^{-k/2}k!I_{n,k}(f_{k})\)^{2M}\le\sum_{\gamma\in\Gamma(k,M)}
n^{-|W(\gamma)|/2}\cdot |F_\gamma|\le\sum_{p=0}^{kM} |\Gamma(k,M,p)|
\eta^{p}2^{p}(kM)^{-p}\sigma^{2M}.
$$
 
Hence by formula (6.1)
$$
\align
E\(n^{-k/2}k!I_{n,k}(f_{k})\)^{2M}&\le
A\(\frac2e\)^{kM}(kM)^{kM}\sigma^{2M}\sum_{p=0}^{kM}\binom{2kM}{2p}
\(2\sqrt{2\eta}\)^{2p}\\
&\le A\(\frac2e\)^{kM}(kM)^{kM}\sigma^{2M}
\(1+2\sqrt{2\eta}\)^{2kM}
\endalign
$$
if $kM\le \eta n\sigma^2$. Thus we have proved Proposition~B with
$C=2\sqrt2$.
 
\medskip\noindent
{\it Proof of Theorem 3.}\/ We can write by the Markov inequality
and Proposition~B with $\eta=\frac{kM}{n\sigma^2}$ that
$$
\aligned
P(k!n^{-k/2}|I_{n,k}(f)|>u)&\le \frac{E\(k!n^{-k/2}I_{n,k}(f)\)^{2M}}
{u^{2M}}\\
&\le A\(\frac1e\cdot 2kM\(1+C\frac{\sqrt{kM}}{\sqrt n\sigma}\)^2
\(\frac\sigma u\)^{2/k}\)^{kM}
\endaligned \tag6.2
$$
for all integers $M\ge M_0$ with some $M_0=M_0(A)$.
 
We shall prove relation (1.8) with the help of estimate (6.2)
first in the case $D\le\frac u\sigma \le n^{k/2}\sigma^k$ with a
sufficiently large constant $D=D(k,C)>0$ depending on $k$ and the
constant $C$ in (6.2). To this end let us introduce the numbers
$\bar M$,
$$
k\bar M=\frac12\(\frac u\sigma\)^{2/k}\frac1
{1+B\frac{ \(\frac u\sigma\)^{1/k}}{\sqrt n\sigma}}
=\frac12\(\frac u\sigma\)^{2/k}\frac1
{1+B\(u n^{-k/2}\sigma^{-(k+1)}\)^{1/k}}
$$
with a sufficiently large number $B=B(C)>0$ and $M=[\bar M]$,
where $[x]$ means the integer part of the number $x$.
 
Observe the $\sqrt{k\bar M}\le\(\frac u\sigma\)^{1/k}$, $\frac{\sqrt
{k\bar M}}{\sqrt n\sigma}\le \(u n^{-k/2}\sigma^{-(k+1)}\)^{1/k}\le 1$,
and
$$
\(1+C\frac{\sqrt{k\bar M}}{\sqrt n\sigma}\)^2\le
1+B\frac{\sqrt{k\bar M}}{\sqrt n\sigma}\le 1+B\(u n^{-k/2}
\sigma^{-(k+1)}\)^{1/k}
$$
with a sufficiently large $B=B(C)>0$ if $\frac u\sigma\le
n^{k/2}\sigma^k$.  Hence
$$
\aligned
\frac1e\cdot 2kM\(1+C\frac{\sqrt{kM}}{\sqrt n\sigma}\)^2
\(\frac\sigma u\)^{2/k}&\le
\frac1e\cdot 2k\bar M\(1+C\frac{\sqrt{k\bar M}}{\sqrt n\sigma}\)^2
\(\frac\sigma u\)^{2/k}  \\
&\le \frac1e\cdot \frac{\(1+C\frac{\sqrt{k\bar M}}{\sqrt n\sigma}\)^2}
{1+B\(u n^{-k/2} \sigma^{-(k+1)}\)^{1/k}}\le\frac1e
\endaligned \tag6.3
$$
if $\frac u\sigma\le n^{k/2}\sigma^k$. If the inequality
$D\le\frac u\sigma$ also holds with a sufficiently large $D=D(B,k)>0$,
then $M\ge M_0$, and the conditions of inequality (6.2) hold. This
inequality together with inequality (6.3) yield that
$$
P(k!n^{-k/2}|I_{n,k}(f)|>u)\le A e^{-kM}\le  Ae^k e^{-k\bar M}
$$
if
$D\le\frac u\sigma \le n^{k/2}\sigma^k$, i.e. inequality (1.8)
holds in this case with a pre-exponential constant $Ae^k$. By
increasing the pre-exponential constant $Ae^k$ in this inequality we
get that relation (1.8) holds for all $0\le\frac u\sigma \le
n^{k/2}\sigma^k$. Thus Theorem~3 is proved.
\medskip
Let us observe that the above calculations show that the constant $B$
in formula (1.8) can be chosen independently of the order $k$ of the
$U$-statistics $I_{n,k}(f)$.
 
\beginsection Appendix. The proof of Theorem B.
 
{\it The proof of Theorem B.}\/ Let us consider all possible sets
$\{(u_1,u'_1),\dots,(u_l,u'_l)\}$, $1\le l\le \min(k_1,k_2)$
containing such pairs of integers for which $u_s\in\{1,\dots,k_1\}$,
$u'_s\in\{1,\dots,k_2\}$, $1\le s\le l$, all points $u_1,\dots,u_l$
are different, and the same relation holds for the points
$u'_1,\dots,u'_l$, too. Let us correspond the diagram
containing two rows $(1,1),\dots,(1,k_1)$ and $(2,1),\dots,(2,k_2)$
and the edges connecting the vertices $(1,u_s)$ and $(2,u'_s)$,
$1\le s\le l$ to the set of pairs $\{(u_1,u'_1),\dots,(u_l,u'_l)\}$,
and let $\bar\Gamma(k_1,k_2)$ denote the set of all (non-coloured)
diagrams we can obtain in such a way. Let us consider the product
$k_1!I_{n,k_1}(f)k_2!I_{n,k_2}(g)$, and rewrite it in the form of
the sum we get by carrying out a term by term multiplication in this
expression. Let us put the terms we get in such a way into disjoint
classes indexed by the elements of the diagrams
$\bar\gamma\in\bar\Gamma(k_1,k_2)$ in the following way : A product
$f(\xi_{j_1},\dots,\xi_{j_{k_1}})g(\xi_{j'_1},\dots,\xi_{j'_{k_1}})$
belongs to the class indexed by the graph
$\bar\gamma\in\bar\Gamma(k_1,k_2)$ with edges
$\{((1,u_1),(2,u'_1)),\dots,((1,u_l),(2,u'_l))\}$ if
$j_{u_s}=j'_{u'_s}$, $1\le s\le l$, for the indices of the random
variables appearing in the above product, and no more coincidence
may exist between the indices $j_1,\dots,j_{k_1},j'_1,\dots,j_{k_2}$.
With such a notation we can write
$$
n^{-(k_1+k_2)/2}k_1!I_{n,k_1}(f)k_2!I_{n,k_2}(g)
=\sumn_{\bar\gamma\in\bar\Gamma\;\;\;}n^{-(k_1+k_2)/2}k(\bar\gamma)!
I_{n,\bar k(\bar\gamma)}(\overline {f\circ g})_{\bar\gamma}), \tag A1
$$
where the functions $(\overline{f\circ g})_{\bar\gamma})$ are
defined in formulas (4.1) and (4.4). (Observe that although formula
(4.4) was defined by means of coloured diagrams, the colours played
no role it. The formula remains meaningful, and does not change if
we replace the coloured diagram $\gamma$ by the diagram $\bar\gamma$
we get by omitting the colours of its edges.) The quantity
$\bar k(\bar\gamma)$ equals the number of such vertices of
$\bar\gamma$ from the first row from which no edge starts plus the
number of vertices in the second row, and the notation
$\sum^{\prime(n)}$ means that summation is taken only for such
diagrams $\bar\gamma\in\bar\Gamma$ for which
$n\ge \bar k(\bar\gamma)$.
 
Let the set $V_1=V_1(\bar\gamma)$ consist of those vertices
$(1,u_1)=(1,u_1)_\gamma$,\dots, $(1,u_{s_1})=(1,u_{s_1})_\gamma$
of the first row $\{(1,1),\dots,(1,k_1)\}$ of the diagram
$\bar\gamma$ from which no edge starts, and let $V_2=V_2(\bar\gamma)$
contain the vertices $(2,v_1)=(2,v_1)_\gamma$,\dots,
$(2,v_{s_1})=(2,v_{s_2})_\gamma$ from the second row
$\{(2,1),\dots,(2,k_2)\}$ of $\gamma$ from which no edges start.
Then $\bar k(\bar\gamma)=s_1+k_2$, and the function
$(\overline {f\circ g})_{\bar\gamma}$ has arguments of the form
$x_{(1,u_p)}$, $(1,u_p)\in V_1$ and $x_{(2,v)}$, $1\le v\le k_2$.
 
Relation (A1) is not appropriate for our goal, since the functions
$(\overline {f\circ g})_{\bar\gamma}$ in it may be non-canonical.
Hence we apply Hoeffding's decomposition for the $U$-statistics
$I_{n,k(\gamma)}(\overline {f\circ g})_{\bar\gamma}$ in formula~(A1)
to get the desired representation for the product of degenerate
$U$-statistics. Actually some special properties of the function
$(\overline {f\circ g})_{\bar\gamma}$ enables us to simplify a
little bit this decomposition.
 
To carry out this procedure let us observe that a function
$f(x_{u_1},\dots,u_{u_k})$ is canonical if and only if
$P_{u_l}f(x_{u_1},\dots,x_{u_k})=0$ with the operator $P_{u_l}$
defined in (4.2) for all indices $u_l$. Beside this, the condition
that the functions $f$ and $g$ are canonical implies the
relations $P_{(1,u)}(\overline {f\circ g})_{\bar\gamma}=0$ for
$(1,u) \in V_1$ and $P_{(2,v)}(\overline {f\circ g})_{\bar\gamma}=0$
for $(2,v)\in V_2$. Moreover, these relations remain valid if we
replace the functions $(\overline {f\circ g})_{\bar\gamma}$ by such
functions which we get by applying the product of some transforms
$P_{(2,v)}$ and $Q_{(2,v)}$, $(2,v)\in\{(2,1),\dots,(2,k_2)\}
\setminus V_2$ for them  with the transforms $P$ and $Q$ defined in
formulas (4.2) and (4.3). (Here we applied such transforms $P$ and
$Q$ which are indexed by those vertices of the second row of
$\bar\gamma$ from which some edge starts.)
 
Beside this, the transforms $P_{(2,v)}$ or $Q_{(2,v)}$ are
exchangeable with the operators $P_{(2,v')}$ or $Q_{(2,v')}$ if
$v\neq v'$, $P_{(2,v)}+Q_{(2,v)}=I$, where $I$ denotes the identity
operator, and $P_{(2,v)}Q_{(2,v)}=0$, since
$P_{(2,v)}Q_{2,v}=P_{(2,v)}-P^2_{(2,v)}=0$. The above relations
enable us to make the following decomposition of the function
$(\overline {f\circ g})_{\bar\gamma}$ to the sum of canonical
functions (just as it is done in the Hoeffding decomposition): Let
us introduce the class of those coloured diagram $\Gamma(\bar\gamma)$
which we can get by colouring all edges of the diagram $\gamma$
either with colour~1 or colour~$-1$. Some calculation shows that
$$
(\overline {f\circ g})_{\bar\gamma}=
\prod_{(2,v)\in\{(2,1),\dots,(2,k_2)\}\setminus V_2}
(P_{(2,v)}+Q_{(2,v)})(\overline {f\circ g})_{\bar\gamma}=
\sum_{\gamma\in\Gamma(\bar\gamma)}(f\circ g)_\gamma, \tag A2
$$
where the function $(f\circ g)_\gamma$ is defined in formula (4.5).
We get the right-hand side of relation (A2) by carrying out the
multiplications for the middle term of this expression, and exploiting
the properties of the operators $P_{(2,v)}$ and $O_{(2,v)}$. Moreover,
these properties also imply that the functions $(f\circ g)_\gamma$
are canonical functions of their variables $x_{(1,u)}$, $(1,u)\in V_1$
and $x_{(2,v)}$, $(2,v)\in B_{(b,-1)}(\gamma)\cup V_2$. (We preserve
the notation of the main part by which $B_{(b,1)}(\gamma)$ and
$B_{(b,-1)}(\gamma)$ denote the sets of those vertices $(2,j)$ of
the second row of the coloured diagram $\gamma$ from which an edge
of colour~1 or colour~$-1$ starts.) Indeed, the above properties of
the operators $P_{(2,v)}$ and $Q_{(2,v)}$ imply that
$P_{(1,u)}(f\circ g)_\gamma=0$ if $(1,u)\in V_1$, and
$P_{(2,v)}(f\circ g)_\gamma=0$ if $(2,v)\in B_{(b,-1)}(\gamma)\cup V_2$.
 
Let $Z(\gamma)$ denote the set of edges of colour~1,
$W(\gamma)$ the set of edges of colour $-1$ in the coloured
diagram $\gamma$, and let $|Z(\gamma)|$ and $W(\gamma)|$ be their
cardinality. Then $(f\circ g)_\gamma$ is a (canonical) function
with $k(\gamma)=k_1+k_2-(|W(\gamma)|+2|Z(\gamma)|)$ variables, and
formula (A2) implies the following representation of the
$U$-statistic $I_{n,\bar k(\bar\gamma)}\(\overline {f\circ
g})_{\bar\gamma}\)$ in the form of a sum of degenerate
$U$-statistics:
$$
n^{-(k_1+k_2)/2}k(\bar\gamma)! I_{n,\bar k(\bar\gamma)}\((\overline
{f\circ g})_{\bar\gamma}\)=n^{-(k_1+k_2)/2}
\sum_{\gamma\in\Gamma(\bar\gamma)}
J_n(\gamma) n^{|Z(\gamma)|}I_{n,k(\gamma)}\((f\circ g)_\gamma\) \tag A3
$$
with $J_n(\gamma)=1$ if $|Z(\gamma)|=0$, and
$$
J_n(\gamma)=\frac{\prodd_{j=1}^{|Z(\gamma)|}
(n-(k_1+k_2)+|W(\gamma)|+|Z(\gamma)|+j)}
{n^{|Z(\gamma)|}}\quad \text{if } \; |Z(\gamma)|>0. \tag A4
$$
The coefficient $J_n(\gamma)n^{|Z(\gamma)|}$ appeared in formula (A3),
since if we apply the decomposition (A2) for all terms
$(\overline {f\circ g})_{\bar\gamma}(\xi_{j_{(1,u)}},\xi_{j_{(2,v)}},
\; (1,u)\in V_1,\,(2,v)\in\{1,\dots k_2\})$ of the $U$-statistic
$k(\bar\gamma)!I_{n,\bar k(\bar\gamma)}\((\overline {f\circ
g})_{\bar\gamma}\)$, then each term
$(f\circ g)_{\gamma}(\xi_{j_{(1,u)}},\xi_{j_{(2,v)}},\;
(1,u)\in V_1,\,(2,v)\in V_2\cup V_1)$ of the $U$-statistic
$I_{n,k(\gamma)}\((f\circ g)_\gamma\)$ appears
$A_n(\gamma)n^{|Z(\gamma)|}$ times. (This is so, because
$k(\gamma)=k_1+k_2-(|W(\gamma)|+2|Z(\gamma)|)$ variables are
fixed in the term $(f\circ g)_{\gamma}$  from the
$\bar k(\gamma)=k_1+k_2-(|W(\gamma)|+|Z(\gamma)|)$ variables in
the term $(\overline {f\circ g})_{\bar\gamma}$, and to get
formula~(A3) from formula~(A2) the indices of the remaining
$|Z(\gamma)|$ variables can be freely chosen from the indices
$1,\dots,n$, with the only restriction that all indices must be
different.
 
Formula (4.6) follows from relations (A1) and (A3). (To see that
we wrote the right power of $n$ in this formula observe that
$n^{-(k_1+k_2)/2}n^{|Z(\gamma)|}=n^{-k(\gamma)/2}n^{-|W(\gamma)|/2}$.)
 
\medskip
To prove inequality (4.7) in the case $|W(\gamma)|=0$ let us estimate
first the value of the function
$(f\circ g)_\gamma^2(x_{(1,u)},x_{(2,v)},\;(u,1)\in V_1,\,(v,2)\in
V_2)$ by means of the Schwarz inequality. We get that
$$
\aligned
&(f\circ g)_\gamma^2(x_{(1,u)},x_{(2,v)},\;(1,u)\in V_1,\; (2,v)\in
V_2)\\
&\qquad \le \int f^2(x_{(1,u)},x_{(2,v)},\;(1,u)\in V_1,\; (2,v)\in
B_{(b,1)}(\gamma))\prod_{(2,v)\in B_{(b,1)}(\gamma)}\mu(\,dx_{(2,v)})
\\ &\qquad\qquad\qquad\int g^2(x_{(2,v)},\;(2,v)\in V_2\cup
B_{(b,1)}(\gamma), )\prod_{(2,v)\in B_{(b,1)}(\gamma)}
\mu(\,dx_{(2,v)})\\
&\qquad=\prod_{(2,v)\in B_{(b,1)}(\gamma)} P_{(2, v)}
f^2(x_{(1,u)},x_{(2,v)},\;(1,u)\in V_1,\; (2,v)\in B_{(b,1)}(\gamma))\\
&\qquad\qquad\qquad\prod_{(2,v)\in B_{(b,1)}(\gamma)} P_{(2,v)}
g^2(x_{(2,v)},\;(2,v)\in V_2\cup B_{(b,1)}(\gamma))
\endaligned \tag A5
$$
with the operators $P$ defined in formula (4.2).
 
Let us observe that the two functions at the right-hand side of~(A5)
are functions of different arguments. The first of them depends on
the arguments $x_{(1,u)}$, $(1,u)\in V_1$, the second one
on the arguments $x_{(2,v)}$, $(2,v)\in V_2$. Beside this, as the
operators $P$ appearing in their definition are contraction in
$L_1$-norm, these functions are bounded in $L_1$ norm by
$\|f\|_2^2$ and $\|g\|_2^2$ respectively. Because of the above
relations we get formula (4.7) by integrating inequality (A5).
 
To prove inequality (4.8) let us introduce, similarly to
formula (4.3), the operators
$$
\tilde Q_{u_j}h(x_{u_1},\dots,x_{u_r})=h(x_{u_1},\dots,x_{u_r})+
\int h(x_{u_1},\dots,x_{u_r})\mu(\,dx_{u_j}), \quad 1\le j\le r,
\tag A6
$$
in the space of functions $h(x_{u_1},\dots,x_{u_r})$ with coordinates
in the space $(X,\Cal X)$. (The indices $u_1,\dots,u_r$ are all
different.) Observe that both the operators $\tilde Q_{u_j}$ and the
operators $P_{u_j}$ defined in (4.2) are positive, i.e. these
operators map a non-negative function to a non-negative function.
Beside this, $Q_{u_j}\le\tilde Q_{u_j}$, and the norms of the
operators  $\frac{\tilde Q_{u_j}}2$ and $P_{u_j}$ are bounded by 1
both in the $L_1(\mu)$, the $L_2(\mu)$ and the supremum norm.
 
Let us define the function
$$
\aligned
&(\widetilde{f\circ g})_\gamma
\(x_{(1,j)},x_{(2,j')},\,j\in\{1,\dots,k_1\}\setminus
B_u(\gamma), \, j'\in \{1,\dots,k_2\}\setminus B_{(b,1)}\) \\
&\qquad=\prod_{(2,j')\in B_{(b,1)}(\gamma)}  P_{(2,j')}
\prod_{(2,j')\in B_{(b,-1)}(\gamma)}  \tilde Q_{(2,j')} \\
&\qquad\qquad\qquad \overline{(f\circ g)}_\gamma
\(x_{(j,1)},x_{(j',2)},\;j\in\{1,\dots,k_1\}\setminus
B_u(\gamma), \, 1\le j'\le k_2\)
\endaligned \tag A7
$$
with the notation of Section~4 in the main part. We have defined
the function $(\widetilde{f\circ g})_\gamma$ with the help of
$\overline{(f\circ g)}_\gamma$ similarly to the definition of
$(f\circ g)_\gamma$ in (4.5), only we have replaced the operators
$Q_{(2,j')}$ by $\tilde Q_{(2,j')}$ in it.
 
We may assume that $\|g\|_2\le\|f\|_2$. We can write because of
the properties of the operators $P_{u_j}$ and $\tilde Q_{u_j}$
listed above and the condition $\sup|f(x_1,\dots,x_k)|\le1$ that
$$
|(f\circ g)_\gamma|\le (\widetilde{|f|\circ |g|})_\gamma\le
(\widetilde{1\circ |g|})_\gamma, \tag A8
$$
where `$\le$' means that the function at the right-hand side is
greater than or equal to  the function at the left-hand side in all
points, and 1 denotes the function which equals identically~1.
Because of relation (A8) it is enough to show that
$$
\aligned
\|(\widetilde{1\circ |g|})_\gamma\|_2&=
\left\|\prod_{(2,j)\in B_{(b,1)}(\gamma)}  P_{(2,j)}
\prod_{(2,j)\in B_{(b,-1)}(\gamma)}  \tilde Q_{(2,j)} \;\;
g(x_{(2,1)},\dots,x_{(2,k_2)})\right\|_2\\
&\le 2^{|W(\gamma)|}\|g\|_2.
\endaligned \tag A9
$$
to prove relation (4.8). But this inequality trivially holds, since
the norm of all operators $P_{(2,j)}$ in formula (A9) is bounded
by~1, the norm of all operators $\tilde Q_{(2,j)}$ is bounded
by~2 in the $L_2(\mu)$ norm, and $|B_{(b,-1)}|=|W(\gamma)|$.
 
\beginsection References.
 
\item{1.)} Arcones, M. A. and Gin\'e, E. (1993) Limit theorems for
$U$-processes. {\it Ann. Probab.} {\bf 21}, 1494--1542
\item{2.)} de la Pe\~na, V. H. and  Montgomery--Smith, S. (1995)
Decoupling inequalities for the tail-probabilities of multivariate
$U$-statistics. {\it Ann. Probab.}, {\bf 23}, 806--816
\item{3.)} Dudley, R. M. (1998)  {\it Uniform Central Limit
Theorems.}\/ Cambridge University Press, Cambridge U.K.
\item{4.)} Dynkin, E. B. and Mandelbaum, A. (1983) Symmetric
statistics, Poisson processes and multiple Wiener integrals. {\it
Annals of Statistics\/} {\bf 11}, 739--745
\item{5.)} Houdr\'e, C. and Reynaud--Bouret, P. (2004): Exponential
inequalities, with constants, for $U$-statistics of order two.
(Preprint)
\item{6.)} It\^o K.: (1951) Multiple Wiener integral. {\it J. Math.
Soc. Japan}\/  {\bf3}.  157--164
\item{7.)} Major, P. (1981) Multiple Wiener--It\^o integrals. {\it
Lecture Notes in Mathematics\/} {\bf 849}, Springer Verlag, Berlin
Heidelberg, New York,
\item{8.)} Major, P. (2004) An estimate about multiple stochastic
integrals with respect to a normalized empirical measure.
submitted to  {\it Studia Scientarum Mathematicarum Hungarica.}
\item{9.)} Major, P. (2004) A multivariate generalization of
Hoeffding's inequality. Submitted to {\it Ann. Probab.}
\item{10.)} Major, P. (2004) On the tail behaviour of multiple
random integrals. (manuscript for a future Lecture Note)
 
\vskip2truecm \noindent
Abbreviated title: Bernstein's inequality.

\bye